\documentclass[]{article}

\usepackage{graphicx}
\usepackage{cancel}
\usepackage{authblk}
\usepackage{amsmath}
\usepackage{amsthm}
\usepackage{amssymb}
\usepackage{bm,upgreek}

\usepackage[sort&compress,numbers]{natbib}
\usepackage[colorlinks=true,linkcolor=black, citecolor=blue, urlcolor=blue]{hyperref}

\begin{document}

\title{Enriched $C^1$ finite elements for crack problems in simplified strain gradient elasticity} 

\author[1,2]{Y.O. Solyaev}
\author[2]{V.N. Dobryanskiy}

\affil[1]{Institute of Applied Mechanics of Russian Academy of Sciences, Moscow, Russia}
\affil[2]{Moscow Aviation Institute, Moscow, Russia}

\setcounter{Maxaffil}{0}
\renewcommand\Affilfont{\itshape\small}

\date{\today}

\maketitle

\begin{abstract}
We present a new type of triangular $C^1$ finite elements developed for the plane strain crack problems within the simplified strain gradient elasticity (SGE). 
The finite element space contains a conventional fifth-degree polynomial interpolation that was originally developed for the plate bending problems and subsequently adopted for SGE. 
The enrichment is performed by adding the near-field analytic SGE solutions for crack problems preserving $C^1$ continuity of interpolation. This allows us an accurate representation of strain and stress fields near the crack tip and also results in the direct calculation of the amplitude factors of SGE asymptotic solution and related value of J-integral (energy release rate). The improved convergence of presented formulation is demonstrated within mode I and mode II problems. Size effects on amplitude factors and J-integral are also evaluated. It is found that amplitude factors of SGE asymptotic solution exhibit a linear dependence on crack size for relatively large cracks.  
\end{abstract}

\section{Introduction}
\label{int}

Strain gradient elasticity (SGE) enables the derivation of regularized solutions for classical singular problems involving cracks, dislocations, concentrated forces, sharp notches, edges, etc. \cite{Mindlin1964,gourgiotis2009plane, aravas2009plane, gourgiotis2010problem, gourgiotis2018concentrated, gutkin1999dislocations, lazar2006dislocations, solyaev2022elastic}. SGE also provides the refined static and dynamic analysis in structured media accounting for size effects and microstructural contributions \cite{rosi2018validity,giorgio2024second,sessa2024implicit}. Application of regularized SGE solutions for the failure analysis of cracked bodies have been proposed and discussed in Refs. \cite{askes2015understanding,vasiliev2019estimation}. It was shown, that SGE solutions allows to describe the transition from short to long crack regimes and related size effects on nominal strength  in quasi-brittle materials \cite{askes2015understanding,vasiliev2019estimation,vasiliev2021new}. Assessments on the values of additional material constants of SGE (the length scale parameters) for different materials and structures were obtained based on the homogenization methods \cite{bacca2013mindlin,ganghoffer2021variational,solyaev2022self,yang2022verification}, atomistic modelling \cite{maranganti2007novel,lazar2022mathematical} and experimental data \cite{vasiliev2021new,vasiliev2021failure,razavi2023length,rezaei2024procedure}.

Since the strain energy functional of SGE depends on the first and the second derivatives of displacements \cite{Mindlin1964}, the numerical simulations within SGE are usually based on the mixed finite element method (FEM) or the conventional FEM with $C^1$-continuous interpolation of displacements. Other kinds of numerical methods, including the boundary element method \cite{Karlis2010}, the meshless and Trefftz-type methods \cite{lurie2006interphase,Wang2019,Solyaev2020}, implicit analysis \cite{Askes2006} and the non-local operator methods \cite{ren2021nonlocal} have been also proposed. 

Mixed FEM within SGE has been developed by using the so-called three-field approximation (for displacement, gradient of displacement and Lagrange multipliers) \cite{shu1999finite}. Modifications, extensions and technical aspects of implementation of mixed FEM in SGE have been presented in Refs.\cite{Amanatidou2002,Zybell2012,Phunpeng2015,reiher2017finite,papanicolopulos2019novel,shekarchizadeh2022benchmark}. Penalty method for mixed FEM within SGE was developed in Refs. \cite{Zervos2008,Zervos2009}. Recently, the variant of mixed FEM with five-field approximation that guarantees the satisfaction of inf-sup (LBB) condition was developed in application to crack problems \cite{Chirkov2024,chirkov2024mixed,chirkov2025mixed}.

$C^1$ finite elements for gradient theories have been developed initially with the third- and fourth-degree interpolation polynomials that provide the continuity of displacement and its first derivatives between the mesh nodes \cite{xia1996crack}. However, it was shown that such interpolation results in a not stable numerical solutions within SGE. Thus, later the fifth-degree polynomial interpolation was successfully applied to solve the governing equations of SGE \cite{Zervos2001,zervos2001modelling}. This interpolation was initially developed for the plate bending problems \cite{dasgupta1990higher} (the so-called Bell triangle with $C^2$ continuity in the mesh nodes) and demonstrated a good accuracy and convergence in different SGE problems, including the cracks problems \cite{akarapu2006numerical,Papanicolopulos2010}. Moreover, it was shown that the computational costs of $C^1$ FEM and mixed FEM in SGE are almost the same, while the former provides the better accuracy\cite{Zervos2009, abali2019computational}. Extension of $C^1$ FEM has been given for various definitions of shape functions \cite{papanicolopulos2013polynomial}, for three-dimensional SGE problems \cite{papanicolopulos2009three} and for the sub-parametric quadrilateral elements \cite{beheshti2017finite}. The second-order homogenization methods was performed with the use of $C^1$ FEM within SGE in Ref.\cite{yvonnet2020computational}.

In the present paper, we developed the modified triangular $C^1$ finite element for the crack problems of SGE. We follow the classical approach, where the  so-called singular elements are placed around the crack tip and embodies the near-field analytical solutions \cite{benzley1974representation,gifford1978stress}. Within SGE we use the corresponding known asymptotic solutions for the plane strain crack problems \cite{gourgiotis2009plane,aravas2009plane,sciarra2013asymptotic} and include them in the set of shape functions of enriched finite elements. The amplitude factors of SGE asymptotic solution become the additional variables of numerical solution. These amplitudes are found simultaneously with the nodal displacements. Moreover, the known relations between the amplitudes and J-integral within SGE \cite{aravas2009plane,sciarra2013asymptotic} allows us to evaluate directly the influence of non-classical size-effects on the energy release rate. Previously, corresponding analysis has been performed based on analytical full-field solutions for crack problems in Ref. \cite{gourgiotis2009plane}.

We show that proposed enriched formulation provide an improved convergence of $C^1$ FEM. Namely, we can use the minimal size of mesh elements of the order $\ell/10$ ($\ell$ is the material length scale parameter) around the crack tip instead of $\ell/100...\ell/1000$ that should be used within the mixed FEM or $C^1$ FEM in SGE \cite {akarapu2006numerical,Papanicolopulos2010,vasiliev2021new}.
We evaluate the amplitude factors of asymptotic solutions for the mode I and II problems and analyse the dependence of these amplitude factors as well as the related values of J-integral on the crack size.

Note that alternative enriched formulation can be also developed by using mixed FEM within SGE, though from our preliminary studies it follows that enriched mixed FEM provides less efficient solutions in comparison to the presented enriched $C^1$ FEM. Similar observation has been obtained previously for the enriched mixed methods within the classical elasticity \cite{heyliger1989stress}.


\section{Simplified strain gradient elasticity theory}
\label{sge}

Consider an isotropic linear elastic body occupying the region $\Omega$ with  boundary $\partial\Omega$ and with the set of edges $\partial\partial\Omega$. The strain energy density within the strain gradient elasticity (Mindlin Form II) can be presented in the following form \cite{Mindlin1964}:
\begin{equation}
\label{w}
\begin{aligned}
	w(\varepsilon_{ij}, \varepsilon_{ij,k}) = 
	\tfrac{1}{2} C_{ijkl} \varepsilon_{ij} \varepsilon_{kl} +
	 &\, \tfrac{1}{2} G_{ijklmn}\varepsilon_{ij,k} \varepsilon_{lm,n}
\end{aligned}
\end{equation}
where $C_{ijkl}=C_{klij}=C_{jikl}=C_{ijlk}$ is the standard tensor of classical elastic moduli and $G_{ijklmn}=G_{lmnijk}=G_{jiklmn}=G_{ijkmln}$ is the sixth-order tensor of gradient moduli; $\varepsilon_{ij} =  \tfrac{1}{2}(u_{i,j}+u_{j,i})$ is an infinitesimal strain tensor; 
$\varepsilon_{ij,k}$ is the strain gradient tensor;
$u_i$ is the displacement vector at a point with coordinates $x_i$; the comma denotes the differentiation with respect to spatial variables and repeated indices imply summation.

The constitutive equations for the Cauchy stress tensor $\tau_{ij}$ and for the third-order double stress tensor $\mu_{ijk}$ are given by:
\begin{equation}
\label{cet}
\begin{aligned}
	\tau_{ij} = \tau_{ji} = \frac{\partial w}{\partial\varepsilon_{ij}} = C_{ijkl} \varepsilon_{kl}
	= \lambda \delta_{ij} \varepsilon_{kk} + 2 \mu \varepsilon_{ij}
\end{aligned}
\end{equation} 
\begin{equation}
\label{cem}
\begin{aligned}
	\mu_{ijk} = \mu_{jik} = \frac{\partial w}{\partial\varepsilon_{ij,k}} = G_{ijklmn} \,\varepsilon_{lm,n} = \ell^2 C_{ijlm} \varepsilon_{lm,k} = \ell^2 \tau_{ij,k}
\end{aligned}
\end{equation}
where $\lambda$, $\mu$ are the L\'ame constants, $\delta_{ij}$ is Kroneker delta, and  we use the constitutive assumption of the simplified SGE that defines the relations between the standard and gradient elastic moduli ($G_{ijklmn} = \ell^2 C_{ijlm}\delta_{kn}$ \cite{altan1992structure,askes2011gradient}) by using single additional length scale parameter $\ell$.

We consider the simplified gradient theory to reduce the amount of programming for the asymptotic solutions, though, the more general theories can be also involved in the developed numerical method. The physical meaning of the length scale parameter $\ell$ \eqref{cem} can be related to the characteristic size of material microstructure \cite{solyaev2022self}. Also, the values of $\ell$ can be correlated to the length scale parameters of fracture mechanics \cite{askes2015understanding,vasiliev2021failure}.

The variation of the strain energy density is given by:
\begin{equation}
\label{dw}
\begin{aligned}
	\delta w= 
	\frac{\partial w}{\partial\varepsilon_{ij}}\delta \varepsilon_{ij} 
	+ \frac{\partial w}{\partial\varepsilon_{ij,k}}\varepsilon_{ij,k} 
	 = \tau_{ij}\delta \varepsilon_{ij} + \mu_{ijk}\delta\varepsilon_{ij,k}
\end{aligned}
\end{equation}

The variation of the total strain energy is therefore:   
\begin{equation}
\label{dW}
\begin{aligned}
	\delta W= \int_\Omega w \,dv 
	 = \int_\Omega (\tau_{ij}\delta \varepsilon_{ij} 
	 		+ \mu_{ijk}\delta\varepsilon_{ij,k})dv
\end{aligned}
\end{equation}

The variation of the total work done by the external forces within SGE is given by \cite{Mindlin1964}:
\begin{equation}
\label{dWe}
\begin{aligned}
	\delta W^{ext} = 
		\int_\Omega b_i\delta u_i \,dv 
		+ \int_{\partial\Omega} \bar t_i\delta u_i \,ds
		+ \int_{\partial\Omega} \bar m_i n_k\delta u_{i,k} \,dv 
		+ \int_{\partial\partial\Omega} \bar s_i\delta u_i \,dl
\end{aligned}
\end{equation}
 where $b_i$ is the body force, $\bar t_i$ is the surface traction, $\bar m_i$ is the surface double traction, $n_i$ is the unit outward normal vector to the boundary $\partial\Omega$,  and $\bar s_i$ is the edge traction.
 
 Considering the linear elastic material, we state that the total strain energy equals to the total work done by the internal forces. Then, the principle of virtual work can be defined by using \eqref{dW} and \eqref{dWe} in the following form:
\begin{equation}
\label{ww}
\begin{aligned}
	\int_\Omega (\tau_{ij}\delta \varepsilon_{ij} + \mu_{ijk}\varepsilon_{ij,k})dv
	&= \int_\Omega b_i\delta u_i \,dv 
		+ \int_{\partial\Omega} \bar t_i\delta u_i \,ds\\[5pt]
		&+ \int_{\partial\Omega} \bar m_i n_k\delta u_{i,k} \,ds 
		+ \int_{\partial\partial\Omega} \bar s_i\,\delta u_i \,dl
\end{aligned}
\end{equation}

Applying the divergence theorem and the Stokes theorem in \eqref{ww}, one can obtain the statement of SGE boundary value problem  \cite{Mindlin1964, dell2015postulations}:
\begin{equation}
\label{bvp}
\begin{cases}
	\sigma_{ij,j} + b_i=0, 
	\qquad &x_i\in\Omega\\[5pt]
	t_i  = \bar t_i, 
	\quad or \quad
	u_i = \bar{u}_i,
	\qquad &x_i\in\partial\Omega\\[5pt]
	m_i  = \bar m_i
	\quad or \quad
	u_{i,j}n_j = \bar{g}_i, 
	\qquad &x_i\in\partial\Omega\\[5pt]
	s_i  = \bar s_i
	\quad or \quad
	u_{i} = \bar{u}^e_i, 
	\qquad &x_i\in\partial\partial\Omega
\end{cases}
\end{equation} 
where $\sigma_{ij} = \tau_{ij} - \mu_{ijk,k}$ is the total stress tensor; $ \bar{u}_i$ and $\bar{g}_i$ are the displacements and normal gradients of displacements that can be prescribed on the body boundary $\partial\Omega$; $ \bar{u}_i^e$ is the displacement that can be prescribed on the body edges $\partial\partial\Omega$; and tractions are defined via stresses and double stresses within SGE in the following form:
\begin{equation}
\label{deft}
	t_i = \sigma_{ij}n_j + \text{D}_{j}(\mu_{ijk}n_k ) 
	+  (\text{D}_{l}n_l) \mu_{ijk}n_jn_k 
\end{equation} 
\begin{equation}
\label{defm}
	m_i = \mu_{ijk}n_jn_k 
\end{equation} 
\begin{equation}
\label{defs}
	s_i = [\mu_{ijk}c_jn_k]
\end{equation} 
where $\text{D}_{i} = (...)_{,i} - n_i(...)_{,k}n_k$ is the surface gradient operator, brackets denote the jump of the enclosed quantities across the edge $\partial\partial\Omega$; $c_j$ is the co-normal vector that is tangent to surface $\partial\Omega$ and normal to edge $\partial\partial\Omega$ \cite{Mindlin1964}.

In this paper, we consider the plain strain problems assuming that the third component of the displacement vector equals to zero and that all field variables depend on the in-plane coordinates only, i.e. $u_1=u(x,y)$, $u_2=v(x,y)$ and $u_3=0$. In this case, only two equilibrium equations remain non-trivial, the surface boundary conditions are reduced to lines, and the edge boundary conditions are reduced to the corner points of the body projection on $xy$-plane.


\section{Asymptotic solution for crack problems}
\label{sol}

In this subsection, we give the representation of analytical asymptotic solutions of SGE that should be used as an enrichment functions in the considered finite element method. Asymptotic solutions for the plane strain crack problems within the simplified SGE were developed in Refs. \cite{gourgiotis2009plane,aravas2009plane} (for discussion, see also \cite{Lazar2015}). These solutions were given for the first two non-zero terms of asymptotic series with behavior $r^1$ and $r^{3/2}$ ($r$ is the distance from the crack tip). The higher order terms of asymptotic series were also derived within SGE in our recent work \cite{solyaev2024higher}. The classical lower order term $r^{1/2}$ is abandoned within SGE.

For the purpose of present study, it is essential to consider only the leading terms $r^{3/2}$ that define the opening of the crack lips and the raise of strain and stress around the crack tip. These terms also define the singular behavior of the second gradient of displacements (and double stress) around the crack tip and the value of J-integral \cite{aravas2009plane,sciarra2013asymptotic}. The terms with $r^1$ behaviour define the constant strain and stress (generalized T-stress field) and zero double stress at the crack tip \cite{gourgiotis2009plane,sciarra2013asymptotic}. Therefore, these terms are already included in the standard polynomial approximation (linear field of displacement) of the finite elements and we do not need to use them as the enrichment functions. Similarly, the next higher order term $r^2$ is also included in the conventional approximation with fifth degree polynomials. The other higher order terms are out of consideration in the present study, though they can be also included in the more general formulation of enriched FEM. It is important to note, that the use of the mentioned terms of asymptotic series ($r^1$, $r^{3/2}$, $r^2$) allows one to fit the analytic SGE solution in the area around the crack tip, including the crack lips and the plane ahead of the crack tip \cite{solyaev2024higher}. 

SGE asymptotic solution for the Cartesian components of displacement vector $\textbf u = u(r,\theta) \textbf e_x + v(r,\theta) \textbf e_y$ can be presented in the following form (see Appendix A):
\begin{equation}
\label{as}
\begin{aligned}
	u = \frac{1}{4\mu}\sum\limits_{n=1}^4 K_n\, Q_{1n},\qquad 
	v = \frac{1}{4\mu}\sum\limits_{n=1}^4 K_n\, Q_{2n}
\end{aligned}
\end{equation}
where 
\begin{equation*}
\begin{aligned}
	Q_{11} &= r^{3/2} \cos\frac{\theta}{2} 
		\left(
		-4 - 2 \eta + 2(1+2\eta) \cos \theta
		\right) \\
	Q_{12} &= r^{3/2}\cos\frac{\theta}{2} 
		\left(
		\frac{7 + 10 \eta}{3} 
		- \frac{14 + 8 \eta}{3} \cos \theta
		- \cos 2\theta
		\right) \\
	Q_{13} &= r^{3/2} \sin\frac{\theta}{2} 
		\left(
		5-4\eta
		+ 4(4+ \eta) \cos \theta
		+ 3\cos 2\theta
		\right)\\
	Q_{14} &= r^{3/2}\sin\frac{\theta}{2} \left(1+2\cos \theta\right)
\end{aligned}
\end{equation*}
\begin{equation*}
\begin{aligned}
	Q_{21} &= r^{3/2} \sin\frac{\theta}{2} 
		\left(
		-4 + 2 \eta - 2(1-2\eta) \cos \theta
		\right) \\
	Q_{22} &= r^{3/2} \sin\frac{\theta}{2} 
		\left(
		\frac{7 - 2 \eta}{3} 
		+ \frac{14 + 8 \eta}{3} \cos \theta
		- \cos 2\theta
		\right) \\
	Q_{23} &= r^{3/2} \cos\frac{\theta}{2} 
		\left(
		-5-12\eta
		+ 4(4+ 3\eta) \cos \theta
		- 3\cos 2\theta
		\right)\\
	Q_{24} &= r^{3/2}\cos\frac{\theta}{2} \left(1-2\cos \theta\right)
\end{aligned}
\end{equation*}
where $K_n$ ($n=1...4$) are the amplitude factors of asymptotic solution; $Q_{in}(r,\theta)$ ($i=1,2$, $n=1...4$) are the functions that define the distribution of asymptotic solution in polar coordinates, $r=\sqrt{x^2+y^2}$ is the distance from the crack tip placed at the point $(0, 0)$ in local Cartesian coordinates with $x$-axis goes beyond the crack tip and $y$-axis goes perpendicular to the crack front; $\theta$ is the angular coordinate estimated from the local $x$-axis in counterclockwise direction; $\eta=3-4\nu$ is Kolosov constant for plain strain; $\nu$ is Poisson's ratio.

Representation \eqref{as} contains SGE solutions for the mode I and mode II crack problems both. Solution for the mode I is defined by constants $K_1$, $K_2$, while the pure mode II is defined by $K_3$, $K_4$.
Note that in contrast to classical elasticity, asymptotic solution of SGE contains two amplitudes for each term in asymptotic series \cite{gourgiotis2009plane,aravas2009plane}. It can be shown that these amplitudes can be related to the so-called  classical and gradient part of general solution for the displacement field (see Appendix A and \cite{solyaev2024higher}). 


The first derivatives of the presented part of SGE near-field solution \eqref{as} behave as $r^{1/2}$. Therefore, these terms define zero strain at the crack tip. The value of maximum (regular) strain at the crack tip is defined by the lower order terms $r^1$ only within SGE. 
The second derivatives of solution \eqref{as} (the strain gradients) are singular and behave as $r^{-1/2}$ around the crack tip. Thus, it is essential to introduce such terms in the set of finite elements shape functions to provide its correspondence to the exact analytical solution. Without these terms, the singularity in the second derivatives necessitates a very dense mesh around the crack tip to achieve the appropriate accuracy of numerical solution.


\section{Finite element formulation}
\label{fem}

Presented formulation of $C^1$ continuous FEM closely follows the method that was originally developed in Refs. \cite{Zervos2001} and used in Refs. \cite{zervos2001modelling,akarapu2006numerical,Papanicolopulos2010}. We propose the enrichment of the elements of this method to improve the convergence of numerical solutions for SGE crack problems. 

In the derivation of weak form of $C^1$ continuous FEM (subsection 4.1), we adopt the standard matrix notation of the finite element method together with extended Voigt notation including the high-grade field variables of SGE. 
In subsection 4.2, we provide a detailed definitions for the standard shape functions used in the considered method, as we utilize a specific subset of these functions for the enrichment procedure outlined in subsection 4.3.


\subsection{Weak form}

The interpolation of the displacement field in each element is defined by:
 
\begin{equation}
\label{u}
\begin{aligned}
		\textbf u = \begin{Bmatrix} u \\ v\end{Bmatrix} = \textbf N\, \hat{\textbf u}
\end{aligned}
\end{equation}
where $\textbf N$ is the matrix of shape functions, and $\hat{\textbf u}$ is the vector of the nodal degrees of freedom. The representations for $\textbf N$ and $\hat{\textbf u}$ will be clarified in the following subsections 4.2, 4.3.

The components of strain and the strain gradient are given by:
\begin{equation}
\label{ek}
\begin{aligned}
		\bm \upvarepsilon = 
		\begin{Bmatrix} 
			\varepsilon_{11} \\ 
			\varepsilon_{22} \\
			2\varepsilon_{12}
		\end{Bmatrix} 
		= \textbf B_1\, \hat{\textbf u},\qquad
		\bm \upkappa = 
		\begin{Bmatrix} 
			\varepsilon_{11,1} \\ 
			\varepsilon_{11,2} \\ 
			\varepsilon_{22,1} \\
			\varepsilon_{22,2} \\
			2\varepsilon_{12,1} \\
			2\varepsilon_{12,2}
		\end{Bmatrix} 
		= \textbf B_2\, \hat{\textbf u},\qquad
\end{aligned}
\end{equation}
where matrices $\textbf B_1$ and $\textbf B_2$ contain the first and the second derivatives of shape functions:
\begin{equation}
\label{b1b2}
\begin{aligned}
		\textbf B_1 &= \textbf L_1 \textbf N, \qquad 
		\textbf B_2 = \textbf L_2 \textbf N\\[5pt]		
		\textbf L_1 &= 
			\begin{pmatrix}
		\tfrac{\partial}{\partial x} & 0 & \tfrac{\partial}{\partial y}\\[5pt]
		0 & \tfrac{\partial}{\partial y} & \tfrac{\partial}{\partial x}
			\end{pmatrix}^\text T,\\[5pt]	
		\textbf L_2 &= 
			\begin{pmatrix}
	\tfrac{\partial^2}{\partial x^2} & 
	\tfrac{\partial^2}{\partial x\partial y} &
	0&0&
	\tfrac{\partial^2}{\partial y\partial x}&
	\tfrac{\partial^2}{\partial y^2}\\[5pt]
	0&0&
	\tfrac{\partial^2}{\partial y\partial x} & 
	\tfrac{\partial^2}{\partial y^2} &
	\tfrac{\partial^2}{\partial x^2} &
	\tfrac{\partial^2}{\partial x\partial y}			
	\end{pmatrix}^\text T
\end{aligned}
\end{equation}

The constitutive equations for stress \eqref{cet} and double stress \eqref{cem} can be rewritten in the following form:
\begin{equation}
\label{ek}
\begin{aligned}
		\bm \uptau = 
		\begin{Bmatrix} 
			\tau_{11} \\ 
			\tau_{22} \\
			\tau_{12}
		\end{Bmatrix} 
		= \textbf C \,\bm \upvarepsilon  
		= \textbf C \,\textbf B_1\, \hat{\textbf u},\qquad
		\bm \upmu = 
		\begin{Bmatrix} 
			\mu_{111} \\ 
			\mu_{112} \\ 
			\mu_{221} \\
			\mu_{222} \\
			\mu_{121} \\
			\mu_{122}
		\end{Bmatrix} 
		= \textbf A \,\bm \upkappa 
		= \textbf A \,\textbf B_2\, \hat{\textbf u}
\end{aligned}
\end{equation}
where the matrices of material constants can be presented within the simplified SGE as follows:
\begin{equation}
\label{ca}
\begin{aligned}	
		\textbf C &= 
			\begin{pmatrix}
		\lambda+2\mu & \lambda & 0\\
		\lambda & \lambda+2\mu & 0\\
		0 & 0 & \mu
			\end{pmatrix},\\[5pt]	
		\textbf A &= 
			\ell^2\begin{pmatrix}
				\lambda+2\mu & 0 & \lambda & 0 & 0 & 0\\
				0 & \lambda+2\mu & 0 & \lambda & 0 & 0\\
				\lambda & 0& \lambda+2\mu & 0 & 0 & 0\\
				0 & \lambda & 0& \lambda+2\mu & 0 & 0 \\
				0&0&0&0&\mu &0\\
				0&0&0&0&0&\mu
			\end{pmatrix}
\end{aligned}
\end{equation}
and the more general SGE constitutive equations can be defined by using appropriate structure of matrix $\textbf A$. 

By using matrix notation \eqref{u}-\eqref{ca}, we can present the virtual work equation \eqref{ww} in the following form:

\begin{equation}
\label{wwf}
\begin{aligned}
	\left(\int_\Omega (
	\textbf B_1^\text T \textbf C \textbf B_1
	+ \textbf B_2^\text T \textbf A \textbf B_2
	)dv\right) \hat{\textbf u}
	&= \int_\Omega \textbf N^\text T \textbf b \,dv 
		+ \int_{\partial\Omega} \textbf N^\text T \textbf t \,ds\\[5pt]
		&+ \int_{\partial\Omega} \bar{\textbf N}^\text T \textbf m \,ds 
		+ \int_{\partial\partial\Omega} \textbf N^\text T \textbf s \,dl
\end{aligned}
\end{equation}
where the variation of nodal variables ($\delta\hat{\textbf u}$) is cancelled, and $\textbf b$, $\textbf t$, $\textbf m$, and $\textbf s$ are the vectors of prescribed vectors of body force, surface traction, double traction and edge traction, respectively, and $\bar{\textbf N}$ is the matrix of normal gradients of the shape functions on the body boundary.

In the following examples, we assume the absence of body force, double traction and edge traction. The virtual work equation can be defined then by introducing the stiffness matrix $\textbf K$ and the load vector $\textbf f$ in the following form:
\begin{equation}
\label{kuf}
\begin{aligned}
	\textbf K  \,\hat{\textbf u} &= \textbf f\\[5pt]
	\textbf K = \int_\Omega (
	\textbf B_1^\text T \textbf C \textbf B_1
	+ \textbf B_2^\text T \textbf A \textbf B_2
	)&dv, \qquad
	\textbf f = \int_{\partial\Omega} \textbf N^\text T \textbf t \,ds
\end{aligned}
\end{equation}

The residual out-of-balance force is given by:
\begin{equation}
\label{r}
\begin{aligned}
	\textbf r  = \int_\Omega (
	\textbf B_1^\text T \bm\uptau
	+ \textbf B_2^\text T \bm\upmu
	)dv -	\textbf f 
\end{aligned}
\end{equation}
where stresses $\bm\uptau$, $\bm\upmu$ are calculated by using found nodal variables $\hat{\textbf u}$. 


\subsection{Standard $C^1$ finite elements}

The conventional approximation for the displacement field utilizes the Bell triangle elements with fifth degree shape functions. These elements were originally developed within the plate bending problems \cite{dasgupta1990higher} and their efficiency within the plain problems of SGE have been shown \cite{Zervos2001,zervos2001modelling}. 
The elements have 12 degrees of freedom per node, or in total 36 nodal variables, including displacements and all its derivatives of the first and the second order. Explicit representation of the shape functions $N_i$ ($i=1...18$) is given in Appendix B to kept this paper self-contained. 
The ordering of the nodal variables is the following: 

\begin{equation}
\label{ui}
\begin{aligned}
	\hat{\textbf u} = \{
						u_1,\, &u_{1,x},\,  u_{1,y},\,
						u_{1,xx}, \,u_{1,xy}, \,u_{1,yy}, \\
						&v_1 \, ... \,v_{1,yy}, \,
						u_2 \, ... \,v_{2,yy}, \,
						u_3 \, ... \,v_{3,yy}\}^\text T
\end{aligned}
\end{equation}
where we introduce the notation for the nodal values of displacements ($u_i, v_i$, $i=1...3$), displacement gradients ($u_{i,x}, \,u_{i,y}$, $u_{i,x}, \,u_{i,y}$, $i=1...3$) and the second gradients of displacement ($u_{i,xx}, \,u_{i,xy}, \,u_{i,yy}$, $v_{i,xx}, \,v_{i,xy}, \,v_{i,yy}$, $i=1...3$). 

Corresponding structure of the matrix of shape function \eqref{u} is the following:
\begin{equation}
\label{Nm}
\small
\begin{aligned}
	\textbf N = 
	\begin{Bmatrix}
		N_1\,...\,N_6 & 0\,...\,0 &
		N_7\,...\,N_{12} & 0\,...\,0 &
		N_{13}\,...\,N_{18} & 0\,...\,0 \\
		0\,...\,0 & N_1\,...\,N_6 & 
		0\,...\,0 & N_7\,...\,N_{12} & 
		0\,...\,0 & N_{13}\,...\,N_{18}
	\end{Bmatrix}
\end{aligned}
\end{equation}
and it is valid that all functions $N_i$, as well as their first and second derivatives take zero values in nodal points ($(x_i,y_i)$, $i=1...3$) except the following cases:
\begin{equation}
\label{Nxy}
\begin{aligned}
	&(x,y)=(x_1,y_1): \quad
	N_1 = N_{2,x}= N_{3,y} = N_{4,xx} = N_{5,xy} = N_{6,yy} = 1 \\
	&(x,y)=(x_2,y_2): \quad
	N_7 = N_{8,x}= N_{9,y} = N_{10,xx} = N_{11,xy} = N_{12,yy} = 1 \\
	&(x,y)=(x_3,y_3): \quad
	N_{13} = N_{14,x}= N_{15,y} = N_{16,xx} = N_{17,xy} = N_{18,yy} = 1
\end{aligned}
\end{equation}

Note that the shape functions of Bell triangle provide even higher accuracy of interpolation ($C^2$ continuity in the nodes) than it is formally required by the weak form of SGE \eqref{ww}. The inconsistency arises then between the approximation and exact solutions in the crack problems, where the second gradient of displacements are infinite and discontinuous at the crack tip \cite{gourgiotis2009plane}. Thus, the proper enrichment is crucial for these elements to obtain the accurate solutions on the coarse meshes.


\subsection{Enriched $C^1$ finite elements}

The enriched finite element approximation is obtained by adding the near-field analytical solution \eqref{as} for the regions near the crack tip to the conventional approximation. In classical methods, the enrichment is performed remaining $C^0$ continuity of approximation by subtracting of nodal values of additional shape function multiplied by the corresponding shape functions \cite{benzley1974representation}. In the present case, we have to preserve $C^1$ continuity of approximation that can be done as follows:
\begin{equation}
\label{ue}
\begin{aligned}
	u &= \sum_{i=1}^{18} \hat u_i N_i  
	+ \sum\limits_{n=1}^4 K_n\, Q^*_{1n}\\
	v &= \sum_{i=1}^{18} \hat v_i N_i  
	+ \sum\limits_{n=1}^4 K_n\, Q^*_{2n}
\end{aligned}
\end{equation}
where the first sums are related to the standard approximation $\textbf u = \textbf N \,\hat{\textbf u}$ \eqref{u}, \eqref{ui}, \eqref{Nm}, in which the summation is performed over the corresponding components of the vector of nodal variables that are defined as $\hat u_i \in  \{(\hat{\textbf u})_i, i = 1...6,\,13...18,\,25...30\}$ and $\hat v_i \in  \{(\hat{\textbf u})_i, i = 7...12,\,19...24,\,31...36\}$. Therefore, $\hat u_i$, $\hat v_i$ ($i=1...18$) denote the nodal values of displacements along $x$- and $y$- axes, respectively, as well as all their first and second derivatives. 

The second sums in \eqref{ue} are related to added asymptotic solution \eqref{as} with amplitude factors $K_n$ and with specially introduced functions $Q^*_{in}$ that are given by:
\begin{equation}
\label{Qx}
\begin{aligned}
	Q^*_{in} &= Q_{in} 
		- \sum\limits_{j=1}^3 Q_{in}^{(j)} N_j^{(0)}
		- \sum\limits_{j=1}^3 \frac{\partial Q_{in}^{(j)}}{\partial x} N_j^{(1)}
		- \sum\limits_{j=1}^3 \frac{\partial Q_{in}^{(j)}}{\partial y} N_j^{(2)} 
\end{aligned}
\end{equation}
where $Q_{in}$ are the standard functions of polar coordinates that define the asymptotic solution \eqref{as} with the position of the crack tip in the one of the element nodes, i.e. $x_0=x_i$, $y_0=y_i$ ($i=1,2$ or 3 according to the geometry of finite element model and ordering of the nodes). The notations $Q_{in}^{(j)}$, $\frac{\partial Q_{in}^{(j)}}{\partial x}$, $\frac{\partial Q_{in}^{(j)}}{\partial y}$ are introduced the values of functions $Q_{in}$ and their first derivatives at $j$-th node $(x_j,y_j)$ ($j=1...3$). These nodal values are multiplied with the corresponding shape functions of Bell triangle that are collected in the following subsets:
\begin{equation}
\label{Nii}
\begin{aligned}
	&N_j^{(0)}:\qquad
	N_1^{(0)} =  N_1, \quad 
	N_2^{(0)} =  N_7, \quad
	N_3^{(0)} =  N_{13}, \\
	&N_j^{(1)}:\qquad
	N_1^{(1)} =  N_2, \quad 
	N_2^{(1)} =  N_8, \quad
	N_3^{(1)} =  N_{14}, \\
	&N_j^{(2)}:\qquad
	N_1^{(2)} =  N_4, \quad 
	N_2^{(2)} =  N_{10}, \quad
	N_3^{(2)} =  N_{16}, \\
\end{aligned}
\end{equation}

Definition \eqref{Qx}, \eqref{Nii} allows us to preserve the $C^1$ continuity of approximation in the nodes of triangular element. Note that in classical enrichment methods the similar definition is usually introduced in the reduced form $Q^*_{in} = Q_{in} - \sum\limits_{j=1}^3 Q_{in}^{(j)} N_j^{(0)}$, that preserve the continuity of displacements only \cite{benzley1974representation,gifford1978stress}. In the case of SGE, we provide additional re-normalization of functions $Q_{in}$ by subtracting the values of their first derivatives multiplied with appropriate shape functions. The subsets of these shape functions \eqref{Nii} are chosen according to their behavior inside the element (see \eqref{Nxy}):
\begin{itemize}
	\item Functions in subset $N_j^{(0)}$ have unit value at $j$-th node, while they have zero values in the other nodes as well as all their first and second derivatives have zero values in all nodes. 
	\item First derivative along $x$-axis of functions in the subset $N_j^{(1)}$ have unit value at $j$-th node, while these functions and all other their first and all second derivatives have zero values in all nodes.
	\item First derivative along $y$-axis of functions in the subset $N_j^{(2)}$ have unit value at $j$-th node, while these functions and all other their first and all second derivatives have zero values in all nodes. 
\end{itemize}

As a result, the functions $Q_{in}^*$ \eqref{Qx} and all their first derivatives have zero values at the nodes of triangular element. These modified functions provide the enrichment of conventional approximation and its first gradient only inside the element. These functions do not change the meaning of conventional nodal degrees of freedom $u_i$, $u_{i,x}$, $u_{i,y}$, $v_i$, $v_{i,x}$, $v_{i,y}$ ($i=1...3$) in vector \eqref{ui}. In such a way, we preserve the $C^1$ continuity of approximation. At the same time, we exclude the redundant $C^2$ continuity of approximation around the crack tip and make it in line with exact solution that requires discontinuous and infinite second gradient of displacements. Note that this discontinuity in the second derivatives is not abandoned by weak form of SGE \eqref{ww}, and its integrability is provided by the correct form of asymptotic solution \eqref{as} (see \cite{gourgiotis2009plane,aravas2009plane}). 

To evaluate the derivatives of $Q_{in}$ in \eqref{Qx} we can use the standard relations between the polar and Cartesian coordinates, and then the relations between the Cartesian and areal coordinates \eqref{dL}. Then, the desired derivatives $\partial/\partial x$, $\partial/\partial y$ can be estimated in terms of areal coordinates via useful formulas \eqref{A6}. As a result, all quantities in the enriched approximation \eqref{ue} will be given in terms of areal coordinates.

The matrix notation \eqref{ui}, \eqref{Nm} can be extended for the elements with enriched approximation \eqref{ue} as follows. The amplitudes of asympototic solution should be included in the vector of nodal variables:
\begin{equation}
\label{uie}
\begin{aligned}
	\hat{\textbf u} = \{
						u_1,\, &u_{1,x},\,  u_{1,y},\,
						u_{1,xx}, \,u_{1,xy}, \,u_{1,yy}, \\
						&v_1 \, ... \,v_{1,yy}, \,
						u_2 \, ... \,v_{2,yy}, \,
						u_3 \, ... \,v_{3,yy},
						K_1,K_2,K_3,K_4\}^\text T 
\end{aligned}
\end{equation}

The matrix $\textbf N$ should be comprised of the standard shape functions and corresponding parts of analytical solution according to \eqref{ue}:
\begin{equation}
\label{Nme}
\footnotesize
\begin{aligned}
	\textbf N = 
	\begin{Bmatrix}
		N_1\,...\,N_6 & 0\,...\,0 &
		N_7\,...\,N_{12} & 0\,...\,0 &
		N_{13}\,...\,N_{18} & 0\,...\,0 & Q_{11}^*...Q_{14}^*\\
		0\,...\,0 & N_1\,...\,N_6 & 
		0\,...\,0 & N_7\,...\,N_{12} & 
		0\,...\,0 & N_{13}\,...\,N_{18} & Q_{21}^*...Q_{24}^*
	\end{Bmatrix}
\end{aligned}
\end{equation}

All other relations of finite element method \eqref{u}-\eqref{kuf} remain the same for the enriched elements.


\subsection{Integration scheme}

\begin{figure}[t!]
\centering
	(a)\includegraphics[width=0.3\linewidth]{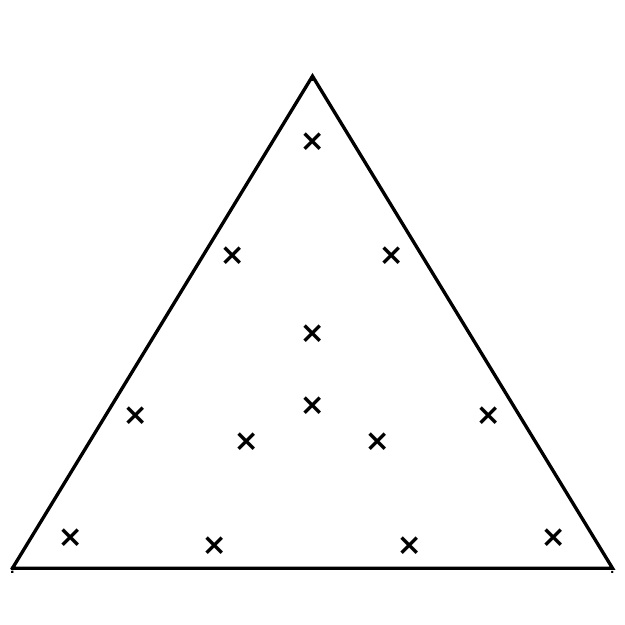}\quad
	(b)\includegraphics[width=0.3\linewidth]{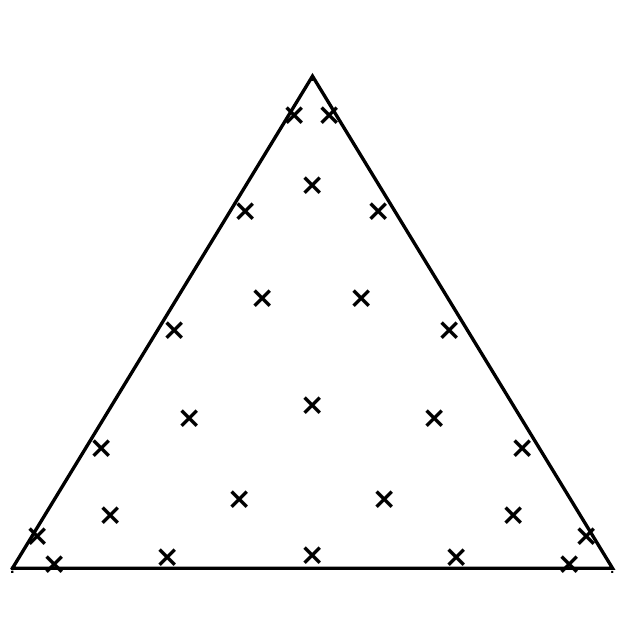}\\
	(c)\includegraphics[width=0.3\linewidth]{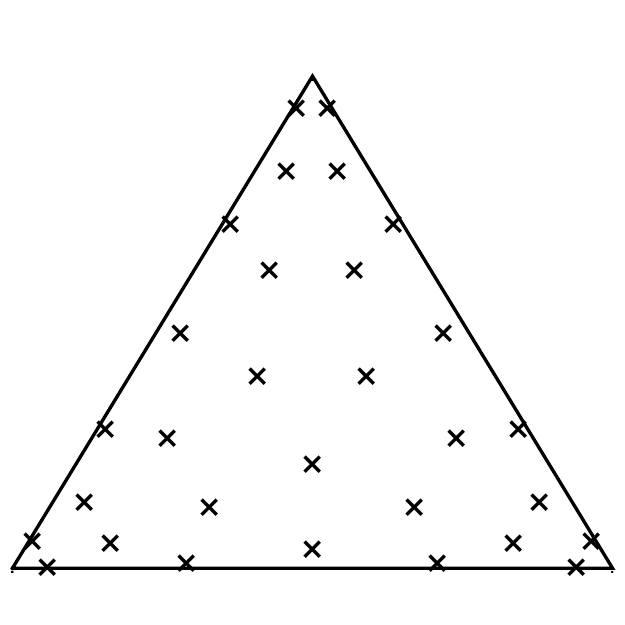}\quad
	(d)\includegraphics[width=0.3\linewidth]{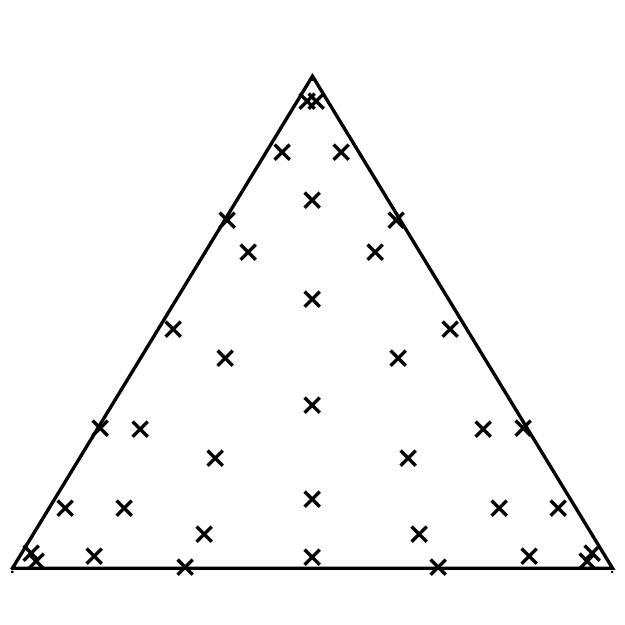}
	\caption{Integration Gauss points used for the conventional elements and enriched elements (a -- reduced 13-points scheme, b -- 25 points, c -- 30 points, b -- 37 points)\label{fig1}}
\end{figure}

The integration for the conventional elements (subsection 4.2) is performed according to the reduced Gauss 13-points scheme (Fig. \ref{fig1}a). The accuracy of this scheme within SGE was shown in Refs. \cite{Zervos2001,Papanicolopulos2010}.  The integration of enriched finite elements is performed by using the full Gauss scheme for the fifth degree polynomials with 25-points (Fig. \ref{fig1}b) and also with the higher-order Gauss scheme with 30 and 37 integration points to evaluate the convergence of algorithm (Fig. \ref{fig1}c, d) \cite{dunavant1985high}. The increased number of integration points can be required by the high first gradients and infinite second gradients of enrichment shape functions around the crack tip that are incorporated in matrices $\textbf B_1$, $\textbf B_2$ in the weak form \eqref{ww}. Integrability of infinite second gradients of enrichment shape functions is provided by the appropriate correct order of asymptotic solution \eqref{as}.


\section{Examples}
\label{fem}

In this section we present the examples of numerical solutions for the mode I and mode II problems for the square region with central crack. We present the comparison between the solutions obtained with the use of conventional and the enriched finite elements. We evaluate the improvement of convergence rate. Also, we evaluate the dependence of amplitude factors of asymptotic solution and related values of J-integral on the length of crack. Presented approach allows us to find these factors together with the displacement solution without additional approaches, developed recently within SGE \cite{chirkov2024computational}.

\begin{figure}[t]
\centering
	(a)\includegraphics[width=0.43\linewidth]{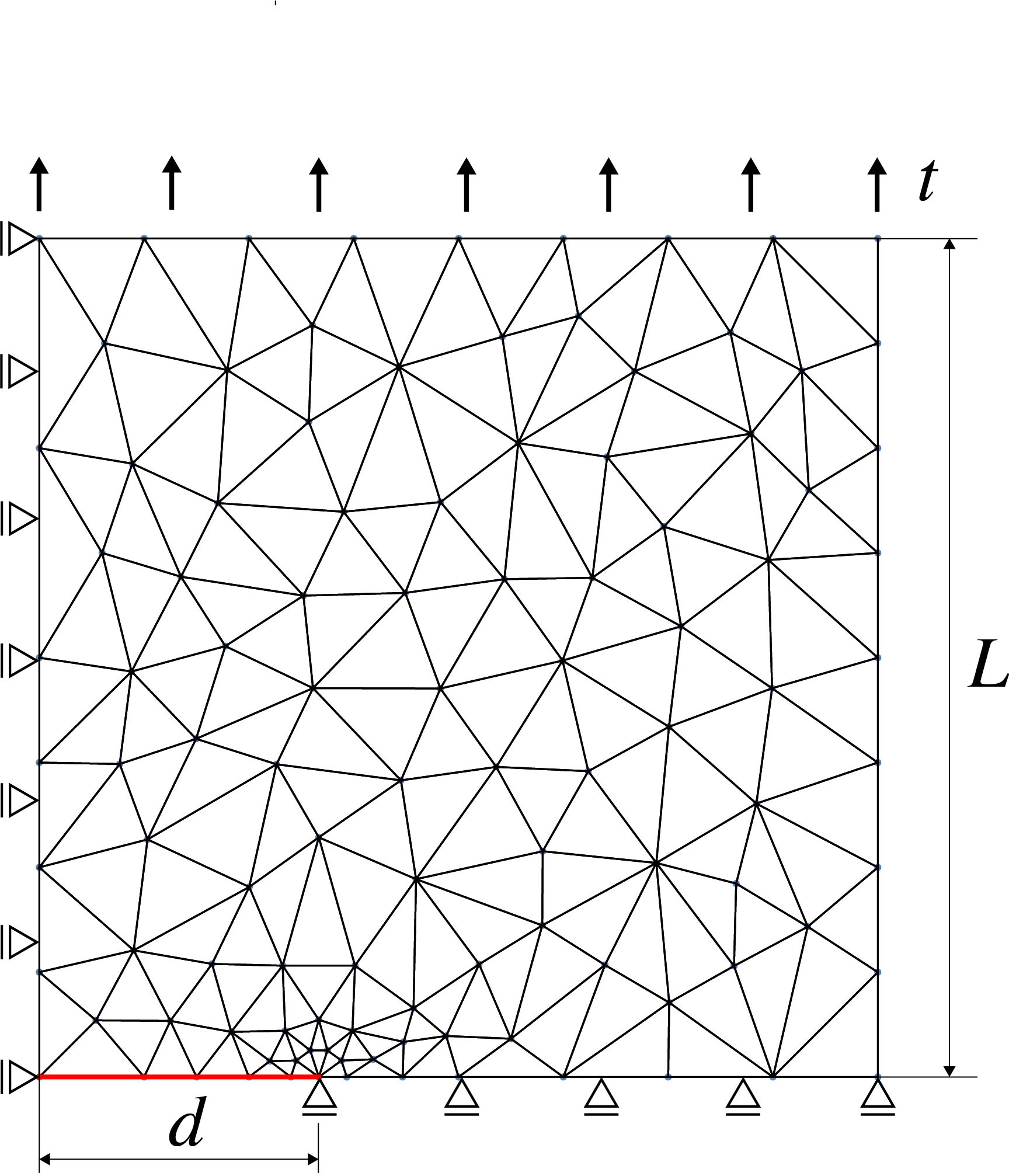}\quad
	(b)\includegraphics[width=0.43\linewidth]{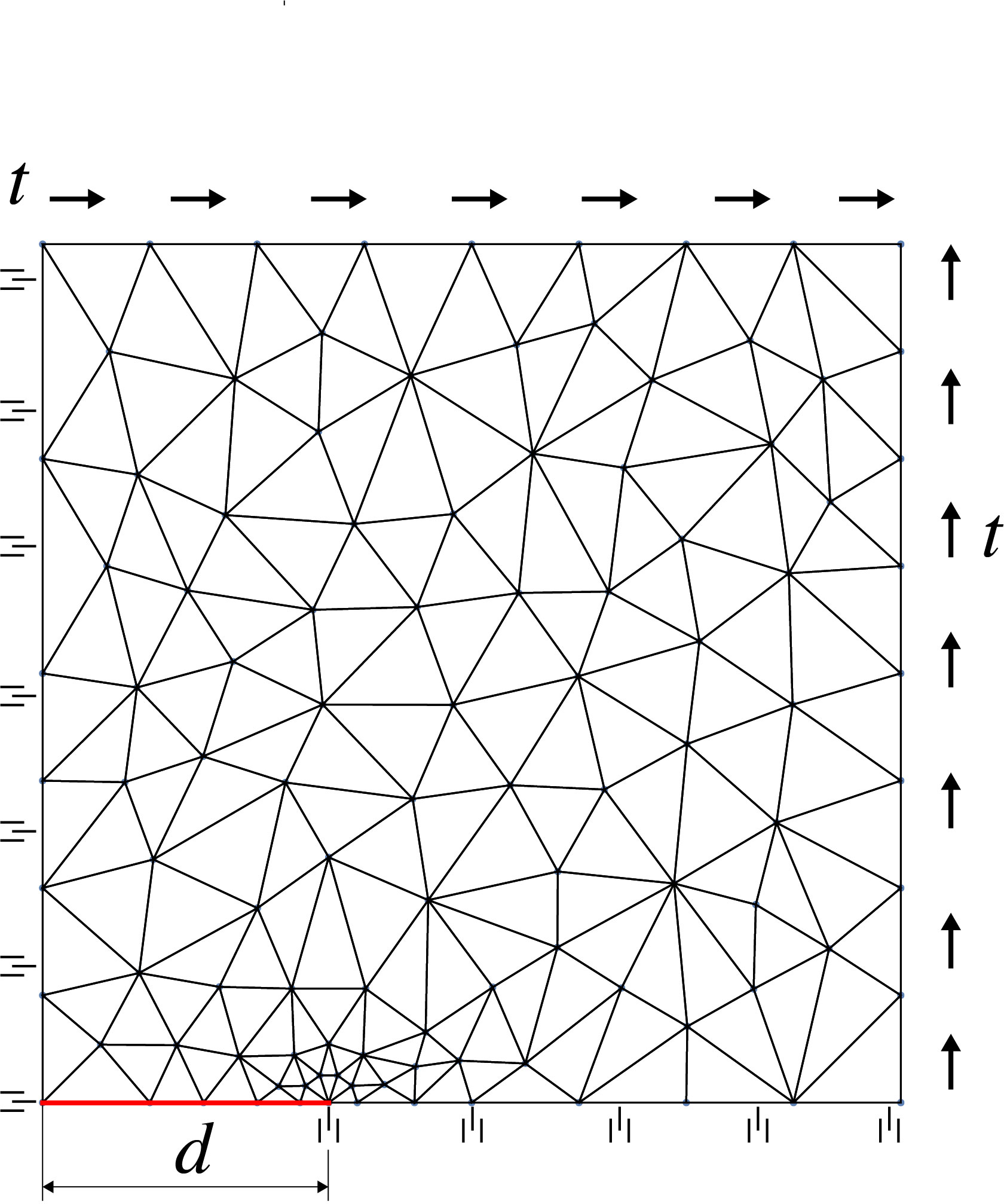}\\[10pt]
	(c)\includegraphics[width=0.45\linewidth]{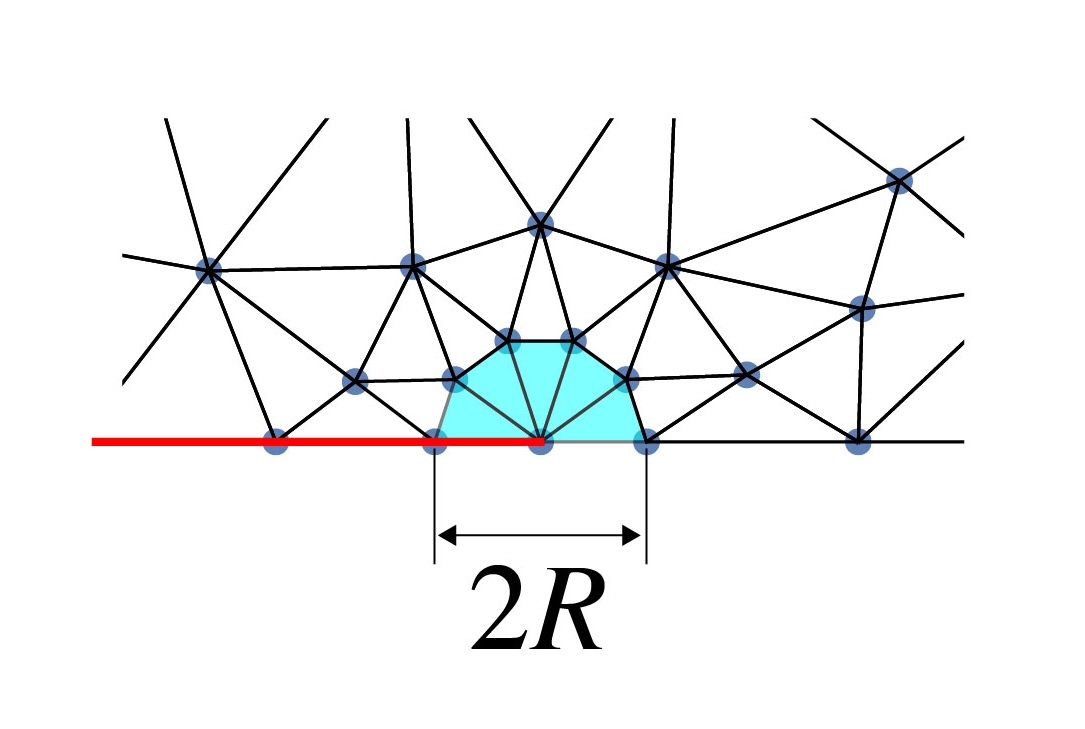}
	\caption{Boundary conditions and examples of finite element mesh used in the test problems.  (a): mode I problem, (c): mode II problem, (d): Enriched elements (blue color) around the crack tip. Cracks are highlighted with red color.}
	\label{fig2}
\end{figure}

The examples of the used mesh and boundary conditions are presented in Fig. \ref{fig2}. We took into account the symmetry of the problems and consider 1/4 of the region. The position of the crack tip is placed at the origin of coordinate system: $x=0,\,y=0$. The crack length is denoted as $d$ and the edge length of the full region is denoted as $2L$. On the planes of symmetry we set the extended form of symmetry/antisymmetry conditions that should be used within $C^1$ FEM for strain gradient theory \cite{Papanicolopulos2010}. Namely, for the mode I problem we prescribe:
\begin{equation}
\label{bcm1}
\begin{aligned}
	x = -d, \,\,0 \leq y\leq L: \quad u =0, \quad u_{,y} =0, \quad  u_{,yy} =0,
			\quad v_{,x} = 0, \quad v_{,xy} =0\\
	 d \leq x\leq L,\,\, y = 0: \quad v=0, \quad v_{,x} =0, \quad  v_{,xx} =0,
			\quad u_{,y} = 0, \quad u_{,xy} =0\\
\end{aligned}
\end{equation}
and for the mode II:
\begin{equation}
\label{bcm1}
\begin{aligned}
	x = -d, \,\,0 \leq y\leq L: \quad v =0, \quad v_{,y} =0, \quad  v_{,yy} =0,
			\quad u_{,x} = 0, \quad u_{,xy} =0\\
	 d \leq x\leq L,\,\, y = 0: \quad u=0, \quad u_{,x} =0, \quad  u_{,xx} =0,
			\quad v_{,y} = 0, \quad v_{,xy} =0\\
\end{aligned}
\end{equation}

Note that the standard definition for the symmetry conditions within SGE implies zero value of normal displacement and zero normal gradient of tangential displacement \cite{lazar2006dislocations,dell2023deformation}. Within the considered numerical method, we also have to prescribe explicitly the absence of tangential derivatives of these quantities since they persist in the set of nodal variables \eqref{ui}, \eqref{uie}. Such definition allows us to fulfil the boundary conditions not only in the mesh nodes, but also between the nodes \cite{Zervos2001,Papanicolopulos2010}. 

In the mode I problem the constant tensile loading $t$ is prescribed at the upper boundary of domain (Fig. \ref{fig2}a). In the mode II problem the tangential loading is prescribed in the correspond positive directions at the upper and right boundaries of the domain (Fig. \ref{fig2}b). The enriched elements are placed around the crack tip and form the circle-like mesh with radius $R$ (Fig. \ref{fig2}c). Due to imposed symmetry conditions, only the one half of this circle persist in the model. The number of enriched elements in this half-circle is denoted as $M$. In the following, we evaluate the influence of $R$ and $M$ on the results of simulations within the considered method. In the example of mesh in Fig. \ref{fig2}c we have $M = 5$.

In the calculations, we used the values of Young's modulus $E = 1$ GPa and Poisson's ratio $\nu = 0.3$ ($\eta = 1.8$). The domain size was $L=1$ m and the load magnitude was $t = 1$ MPa, though all results will be given in dimensionless form for different relative size of crack ($d/L$) and different ratio between the length scale parameter of simplified SGE and the crack length ($\ell/d$). The values of estimated Cauchy stresses will be normalized with respect to the prescribed load level $t$.

\subsection{Comparison between the conventional and enriched elements}

\begin{figure}[t!]
\centering
	(a)\includegraphics[width=0.4\linewidth]{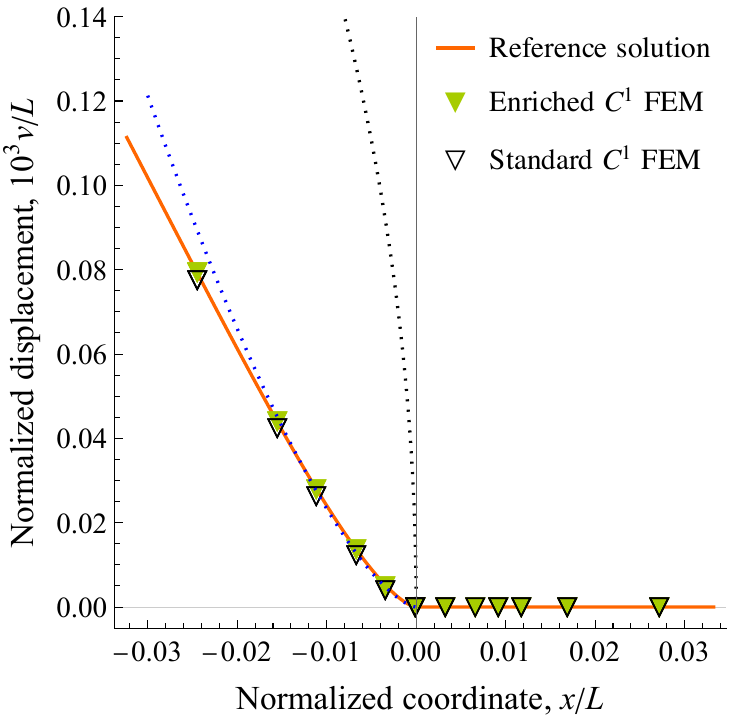}\quad
	(b)\includegraphics[width=0.4\linewidth]{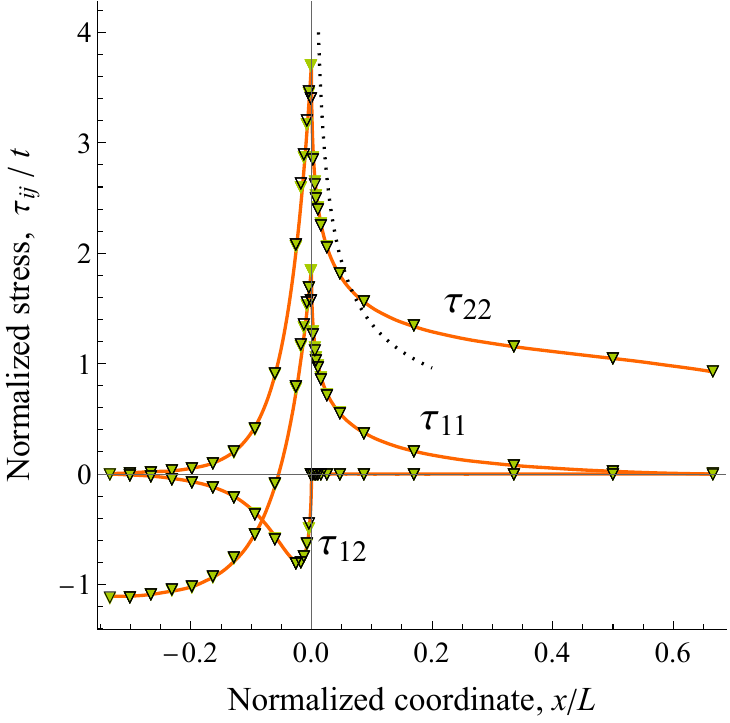}
	(c)\includegraphics[width=0.4\linewidth]{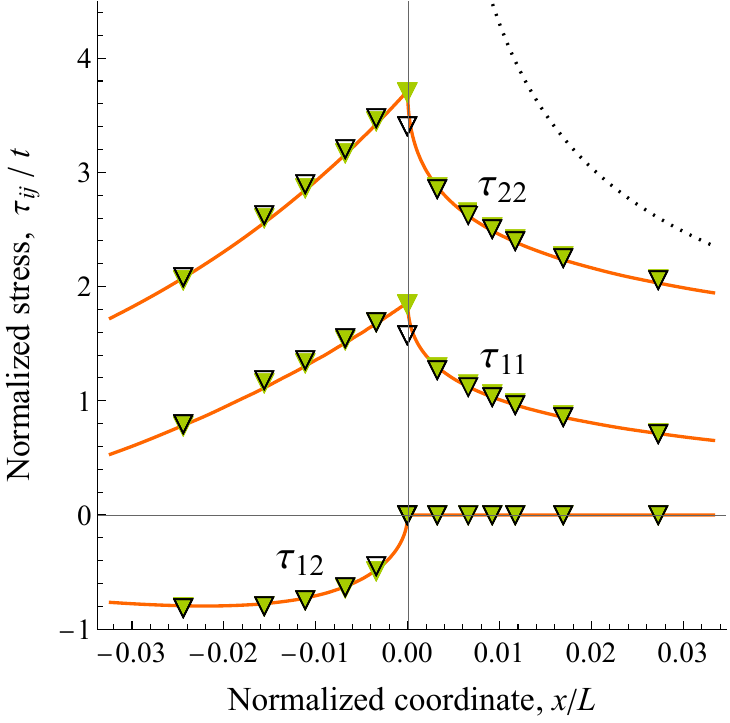}\quad
	(d)\includegraphics[width=0.4\linewidth]{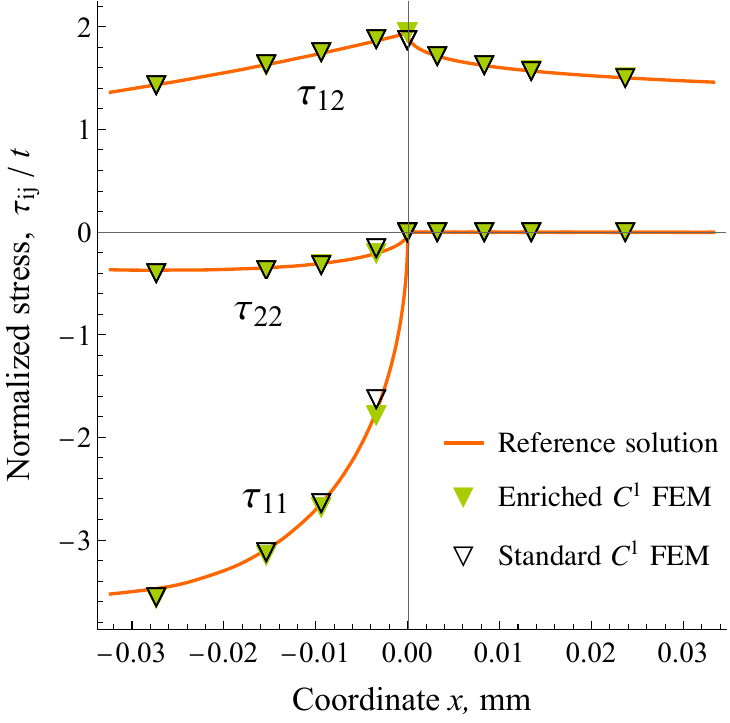}
	\caption{Distribution of field variables along the lower boundary of the domain containing crack. (a) - normal displacements in mode I problem, (b) - Cauchy stresses in mode I problem, (c) - magnified view of figure (b) around the crack tip, (d) - Cauchy stresses in mode II problem. Classical asymptotic solution is shown by black dotted line. In plot (a), SGE asymptotic solution is show by blue dotted line}
	\label{fig3}
\end{figure}

The comparison between the distribution of displacements and Cauchy stresses obtained by using the conventional and the enriched $C^1$ finite elements is presented in Fig. \ref{fig3}. The example is given for the relative crack length $d/L = 3$, the length scale parameter $\ell = d/10$, size of elements around crack $R = \ell/10$ and the number of elements around crack $M=5$. In these figures, we also present the solution that was obtained by using the mixed FEM method implemented in Comsol Multiphysics with the use of the Weak Form PDE interface. We denote this solution as "reference solution" (solid lines in Fig. \ref{fig3}) since it was obtained by using extremely dense mesh with the seed size around the crack tip $R = \ell/1000$.  The implementation technique in Comsol was explained and evaluated in details in Refs. \cite{reiher2017finite,vasiliev2021failure}. The validation of this mixed FEM method based on comparison with different analytical solutions has been also presented \cite{abali2017theory,shekarchizadeh2022benchmark}. We do not compare the present method with some kind of analytical solution since there is no closed form analytical solutions for plane strain crack problems in SGE. The known analytical solution for cracks of finite length has been obtained with the use of numerical integration for the inversion of Fourier transform \cite{gourgiotis2009plane}.

In Fig. \ref{fig3}a we present the crack opening profile around the crack tip for the mode I problem. It can be seen that the cusp-like behavior is realized in SGE solution according to $\sim r^{3/2}$ asymptotic low \eqref{as}. The classical profile ($\sim r^{1/2}$) is presented by black dotted line in Fig. \ref{fig3}a. It can be seen, that all numerical solutions provide close results, though the standard $C^1$ elements provide some underestimation of normal displacements (black markers in \ref{fig3}a). SGE asymptotic solution for the normal displacement (blue dotted line in Fig. \ref{fig3}a) fits well with the full-field numerical solutions. This asymptotic solution is plotted according to the presented form \eqref{as} evaluated at $\theta=\pi$, that is:
$$v =-\frac{|x|^{3/2}}{2\mu}(1+\eta)(K_1 + \tfrac{5}{3}K_2),$$
where the found values of amplitude factors in the considered example are $K_1 = -4.1692$ [MPa$\cdot$m$^{-1/2}$] and $K_2 = -1.3455$ [MPa$\cdot$m$^{-1/2}$]. The dimensions of these amplitudes are different to classical stress intensity factors since the asymptotic behavior of solution is also changed.

In Fig. \ref{fig3}b we present the comparison for Cauchy stress distribution along the lower boundary of the domain ($y=0$) that includes the crack face. Corresponding magnified plot around the crack tip is presented in Fig. \ref{fig3}c. SGE solutions provide regularized stress field (related to the corresponding regularized strain field) with maximum concentration of normal stress $\tau_{22}$. Classical singular asymptotic solution for $\tau_{22}$ is presented by black dotted line in these figures. 
It can be see that exactly at the crack tip, solution obtained with conventional $C^1$ elements has strong deviation from the reference mixed FEM solution and solution obtained with enriched $C^1$ FEM. The reason for this deviation is the continuity of second derivatives of displacement that is imposed by conventional $C^1$ elements and that is violated by the exact SGE solution. Generally say, solution with conventional $C^1$ elements can never reach the exact SGE solution but it tends to this solution with the decrease of the mesh size around the crack tip (see next section). 
Similar result can be observed for the mode II problem solution in Fig. \ref{fig3}d where we show the distribution of stresses in the magnified zone around the crack tip. Here the well observed deviation between the solutions with conventional and enriched elements can be seen not only at the crack tip for the non-zero shear stress $\tau_{12}$ but also for normal stresses that have strong gradients at the crack face side.

Note that in SGE crack problems, Cauchy stresses have non-zero values at the traction-free crack faces ($x<0$ in Fig. 3) since the boundary conditions for traction are prescribed with respect to the combination of Cauchy stresses and gradients of double stresses \eqref{deft}. The non-zero Cauchy stresses at the crack faces can be treated as some kind of cohesion forces that naturally arise in SGE solutions \cite{sciarra2013asymptotic}. These effects were discussed in details in Refs.\cite{aravas2009plane,sciarra2013asymptotic,Papanicolopulos2010,vasiliev2021failure,solyaev2024higher}. The illustrations for the full-field distribution of strains, stresses, and high-order stresses in SGE crack problems can be also found in these works.

\subsection{Convergence analysis}

Convergence analysis is presented for the maximum values of Cauchy stresses estimated at the crack tip within the mode I and II problems. We consider various size of mesh elements around the crack tip $R$, various number of  elements placed around the crack tip $M$ (see Fig. \ref{fig2}) and various order of Gauss integrations scheme for enriched elements (see Fig. \ref{fig1}). Examples of convergence analysis are given for relative crack size $d/L = 1/5$ and for different ratios $d/\ell$, i.e. for different values of the length scale parameter. This parameter can be treated as some characteristic size of material's microstructure so that in the case of small relative values of $d/\ell$ the model corresponds to the strong contribution of gradient effects with lower stress concentration \cite{aravas2009plane,gourgiotis2018concentrated,vasiliev2021failure}. In the case of large values of $d/\ell$ (small values of length scale parameter) the stress concentration increases and it tends to classical singular solution in the limit case when $d/\ell\rightarrow\infty$.

\begin{figure}
\centering
	(a)\includegraphics[width=0.4\linewidth]{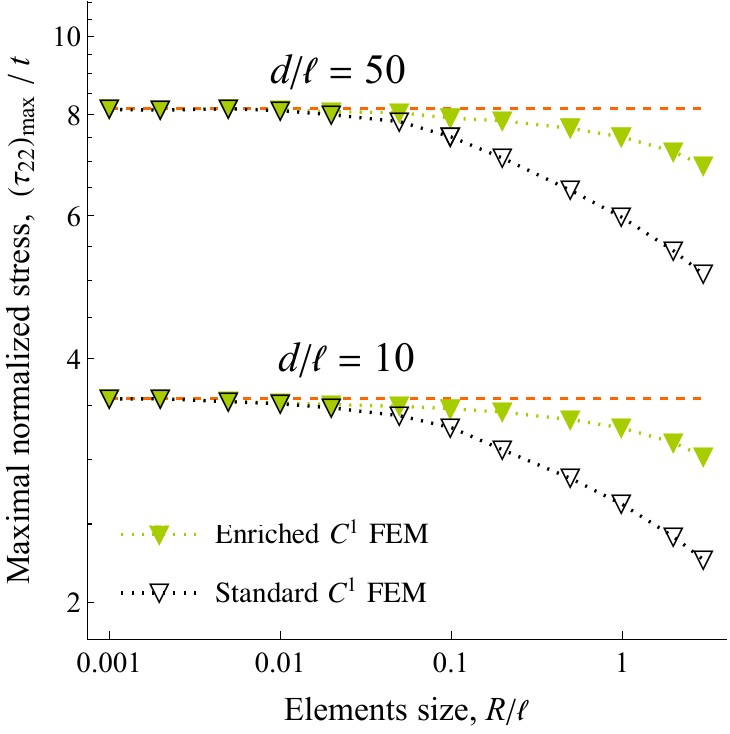}\quad
	(b)\includegraphics[width=0.4\linewidth]{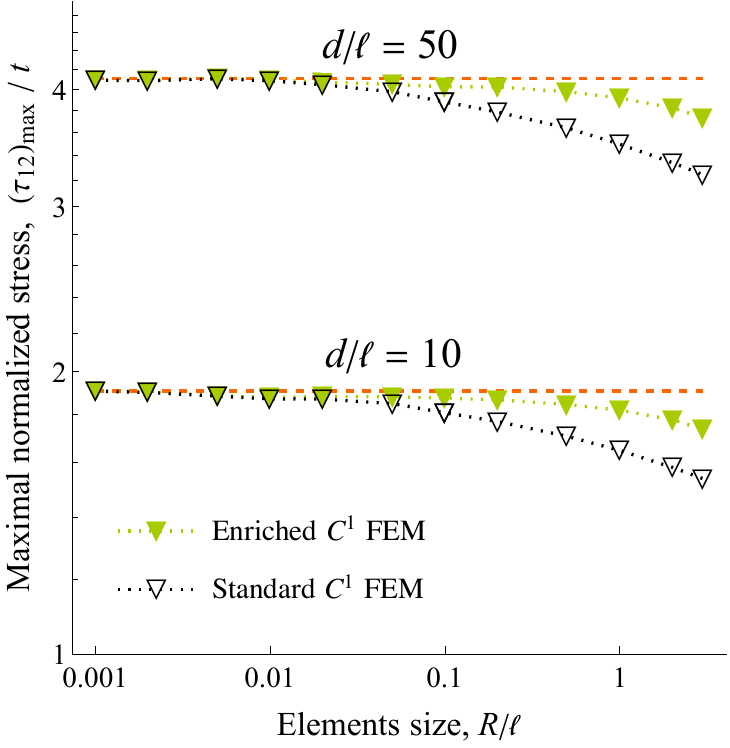}
	\caption{Dependence of maximum stress concentration in mode I (a) and mode II (b) on the relative size of elements placed around the tip of crack. Dashed orange lines correspond to the solutions obtained with the smallest mesh\vspace{2mm}}
	\label{fig4}
\centering
	(a)\includegraphics[width=0.4\linewidth]{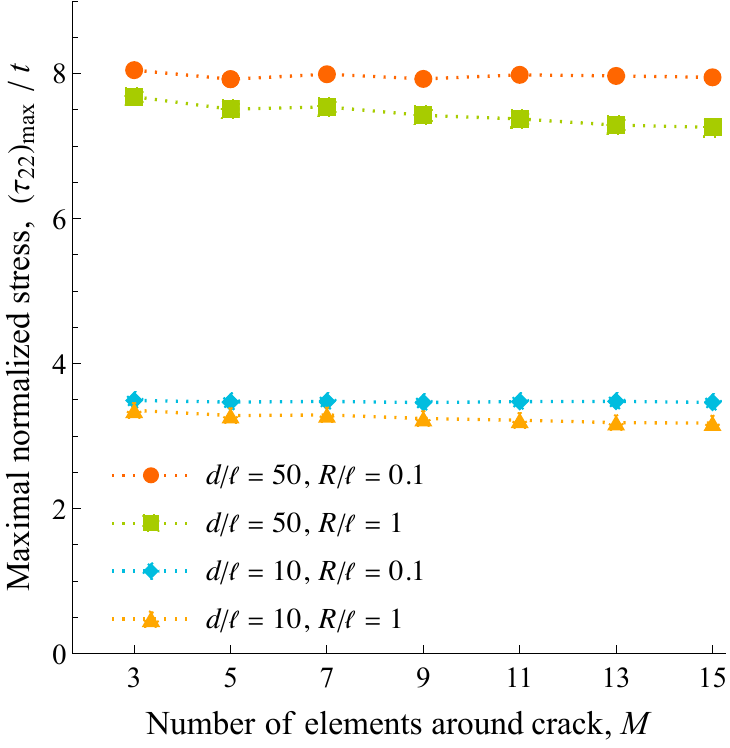}\quad
	(b)\includegraphics[width=0.4\linewidth]{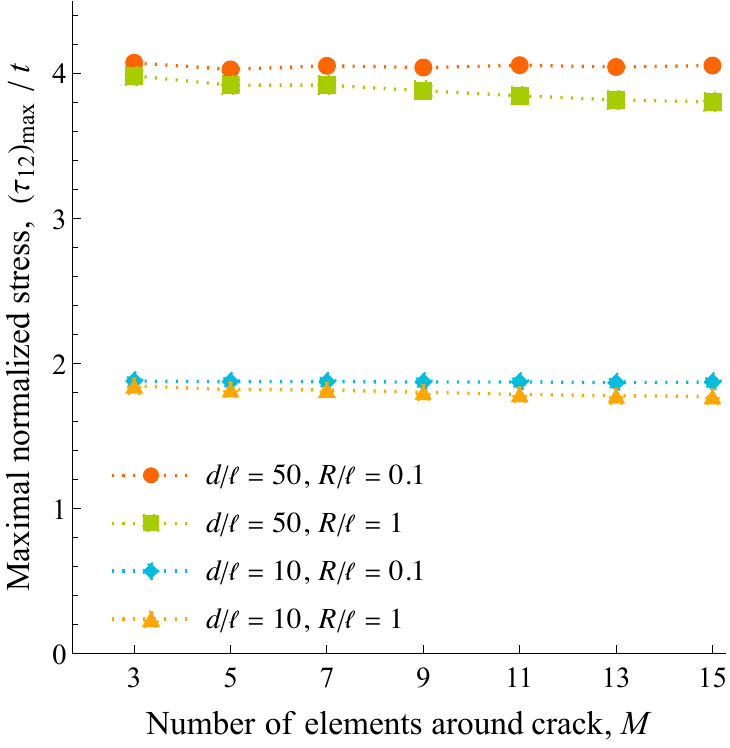}
	\caption{Dependence of maximum stress concentration in mode I (a) and mode II (b) on the number of enriched elements placed around the tip of crack\vspace{2mm}}
	\label{fig5}
\centering
	(a)\includegraphics[width=0.4\linewidth]{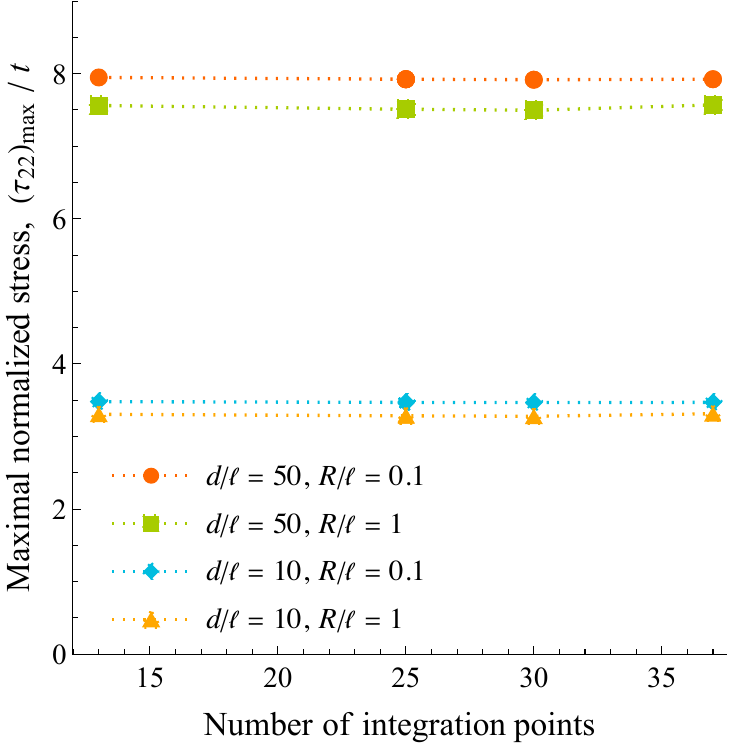}\quad
	(b)\includegraphics[width=0.4\linewidth]{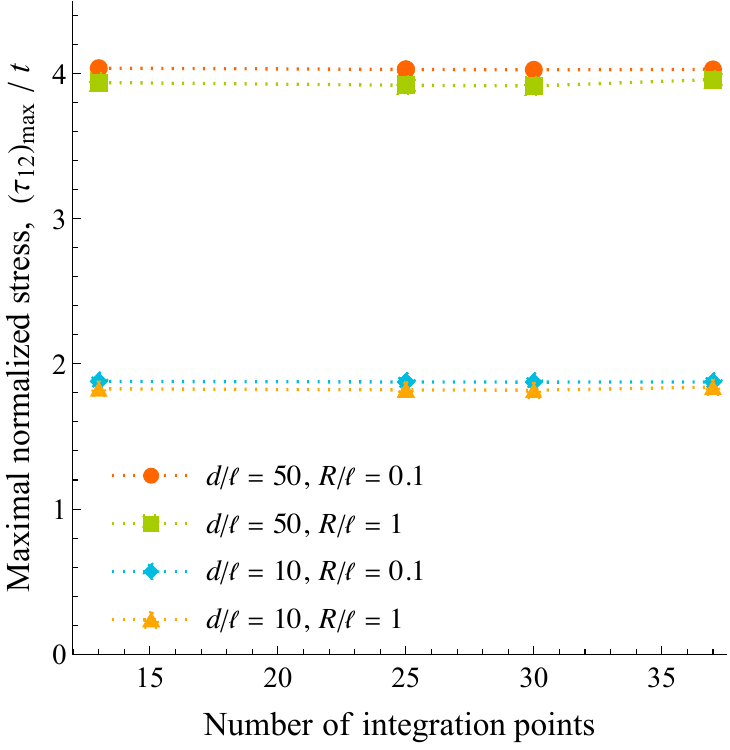}
	\caption{Dependence of maximum stress concentration in mode I (a) and mode II (b) problems on the number of integration points in Gauss quadrature rule  used for the enriched elements}
	\label{fig6}
\end{figure}

Dependence of normalized Cauchy stresses on the relative element size $R/\ell$ is presented in Fig. \ref{fig4}. It can be seen that both solutions with conventional and enriched $C^1$ elements tend to the same asymptotic value with the decrease of element size. However, enriched elements provide better convergence. The error of 2\% can be obtained when using the enriched elements with relative size $R/\ell \approx 0.1$. This error is estimated with respect to the solution obtained by using extremely small size of elements $R/\ell = 0.001$, shown by orange dashed lines in Fig. \ref{fig4}.  To obtain the same accuracy in the model with standard elements, their relative size should be in order smaller, i.e. $R/\ell < 0.01$. The results of these calculation in Fig. \ref{fig4} were obtained by using $M=5$.

The influence of the number of elements placed around the crack tip ($M$) on the maximum stress concentration is shown in Fig. \ref{fig5}. These results are presented for the model with enriched $C^1$ elements. 
It can be seen, that for the chosen enough small relative size of elements ($R/\ell=0.1$) there is no significant dependence of the solution on $M$ (see red an blue lines in Fig. \ref{fig5}). In contrast, for the larger size of enriched elements ($R/\ell=1$, green and orange lines in Fig. \ref{fig5}) there arise an underestimation of stress concentration that increases with increase of elements number $M$. These results are valid for the mode I and II problems both. 

The influence of the order of Gauss integration scheme is presented in Fig. \ref{fig6}. We found that for the presented examples, it is enough to use the reduced scheme with 13 integration point suggested initially for the standard $C^1$ elements (Fig. \ref{fig1}a). The increase of integration points number does not change the results of calculations for different values of the length scale parameter and different size of elements. The number of elements around the crack tip in the examples in Fig. \ref{fig6} was $M=5$. 

Thus, based on the performed analysis we can conclude that it is enough to use rather small size of enriched elements with $R/l<0.1$ and reduced integration scheme to obtain the appropriate accuracy of developed method. The number of elements placed around the crack tip $M>3$ provide stable results of calculations. In the next section, we use these settings ($R/l=0.1$, $M=5$) to evaluate the effects of crack size on the amplitude factors of SGE asymptotic solution.

\subsection{Amplitude factors of asymptotic solution and J-integral}

It is of interest to evaluate the values of amplitude factors of asymptotic solution \eqref{as}. Up to date, there is no systematic data on the values of these factors within SGE. To the best of authors knowledge, assessments were obtained only within the mode I problem in Refs. \cite{aravas2009plane,solyaev2024higher} based on the fitting of asymptotic solution to the full-field numerical solutions. At the same time, the amplitude factors define the value of J-integral (energy release rate). For the mode I and mode II problems the following relations are valid within SGE:
\begin{equation}
\label{J12}
\begin{aligned}
	J_I &= \frac{1+\eta}{8\mu} \pi\ell^2 \left((3 K_1 + K_2)^2 + 8 K_2^2 (\eta+2)\right),\\
	J_{II} &= \frac{1+\eta}{8\mu} \pi\ell^2 \left(72 K_3^2 (\eta+2) 
	+ \frac{9K_4^2}{4 (\eta^2-1)}\right)
\end{aligned}
\end{equation}

These relations can be obtained by using asymptotic solution \eqref{as} (or its representation in polar coordinates \eqref{A1}-\eqref{A6}) and  estimation of the corresponding generalized path independent integral that takes into account the contribution of strain gradient effects (see \cite{aravas2009plane,sciarra2013asymptotic}). The presented result for J-integral \eqref{J12} coincides with those one given in Refs. \cite{aravas2009plane,sciarra2013asymptotic} up to re-normalization of amplitude factors. 

\begin{figure}[b!]
\centering
	(a)\includegraphics[width=0.42\linewidth]{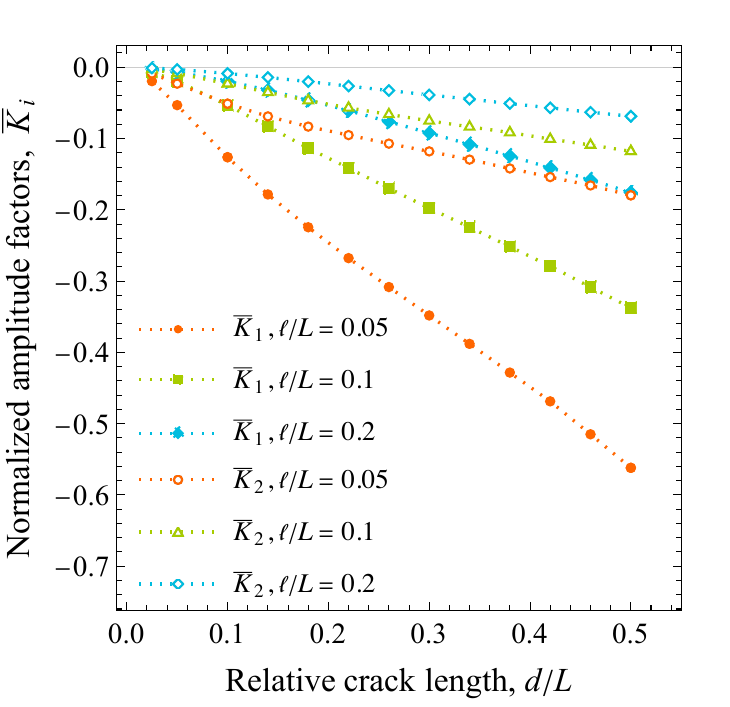}\quad
	   \includegraphics[width=0.4\linewidth]{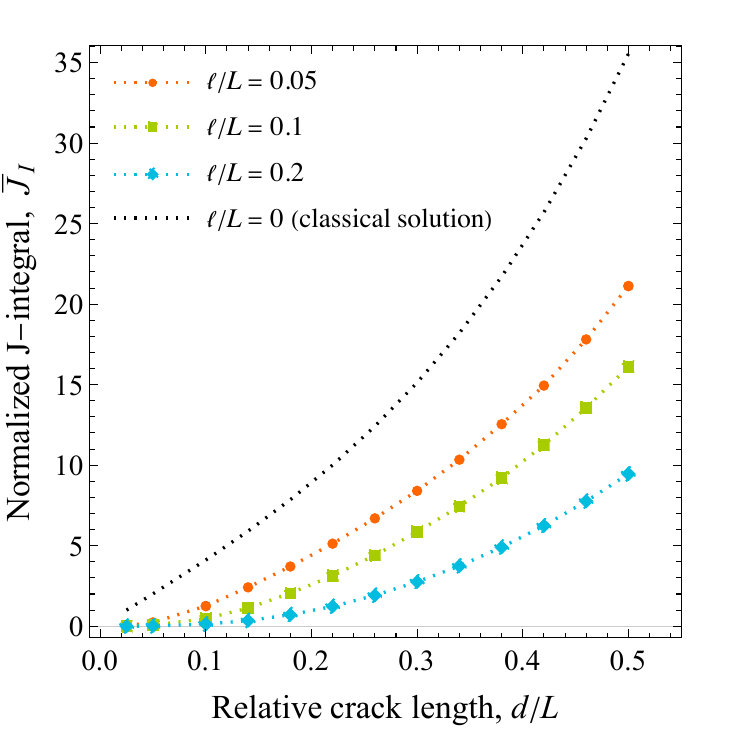}\\
	(b)\includegraphics[width=0.42\linewidth]{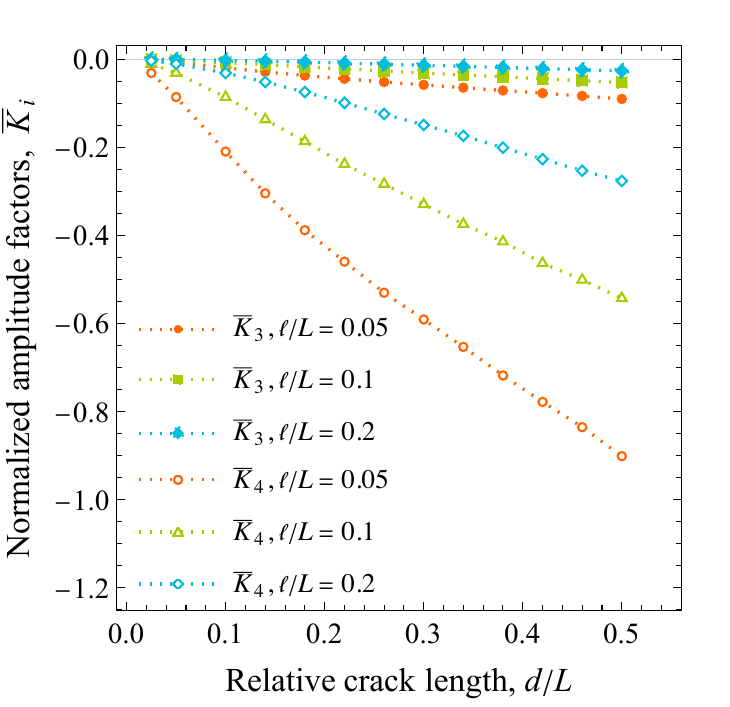}\quad
	   \includegraphics[width=0.4\linewidth]{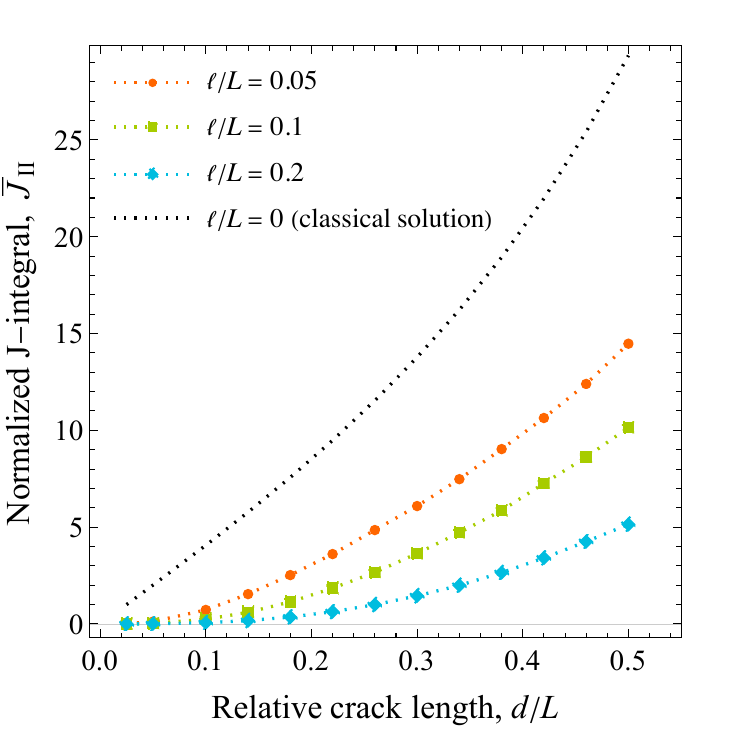}
	\caption{Dependence of amplitude factors of SGE asymptotic solution (left) and J-integral on the relative length of crack under mode I (a) and mode II (b) loading conditions. Classical solutions for J-integral are shown by black dotted line}
	\label{fig7}
\end{figure}

In the examples of calculations, we found the values of amplitude factors $K_i$ and J-integral for the problems with finite size crack of different relative length $d/L$ under mode I and mode II loading conditions. The calculations were performed for different values of the length scale parameter defined with respect to the size of the domain ($\ell/L$). 
In Fig. \ref{fig7} we show the found dependence of amplitude factors and J-integral on the relative length of crack. The values of $K_i$ are normalized on these plots with respect to the dimensional group $t/\sqrt l$ ($t$ is the prescribed magnitude of the load, $l$ is the value of the length scale parameter noted in the legends in Fig. \ref{fig7}). The values of J-integral are normalized with respect to the classical value $J_0$ found for given problems with considered minimal length of crack as follows:
$$\bar J = J/J_0, \quad
J_0 = K_0 \tfrac{1+\eta}{8\mu}
= K_0 \tfrac{1-\nu^2}{E}, \quad K_0 = t \sqrt{\pi d_{min}}, 
\quad d_{min} = L/40$$

It can be seen that the amplitude factors in Fig. \ref{fig7} exhibit the linear dependence on $d/L$ for relatively long cracks in mode I and mode II problems both. The values of $K_i$ are always negative since they describe the rate of stress decrease around the crack tip (the peak value of stress is defined by the lower order linear terms that behave as $\sim r^1$ \cite{aravas2009plane,solyaev2024higher}). The related values of J-integral are found by using Eq. \eqref{J12} and they are compared to the classical J-integral values presented by black dotted lines in Fig. \ref{fig7}. Note that the J-integral in SGE solution is always lower than classical one due to strengthening effects \cite{aravas2009plane,gourgiotis2009plane}. For larger values of the length scale parameter, strengthening effects become more pronounced and the strain energy release rate reduces (see blue lines in Fig. \ref{fig7}). Thus, the presented method with enriched finite elements can be useful for evaluation of energy-based fracture criteria and corresponding experimental validation of SGE theory for quasi-brittle materials exhibiting size effects in fracture. The data on the dependence of critical J-integral on the size of crack and the size of specimen can be used for identification of the length scale parameter as an additional material constant. 

\begin{figure}[b!]
\centering
	(a)\includegraphics[width=0.4\linewidth]{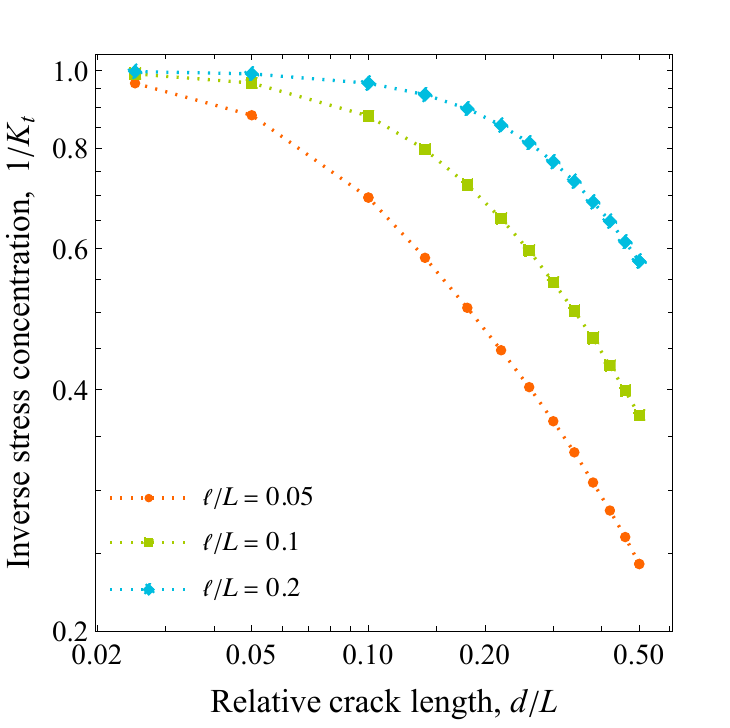}\quad
	(b)\includegraphics[width=0.4\linewidth]{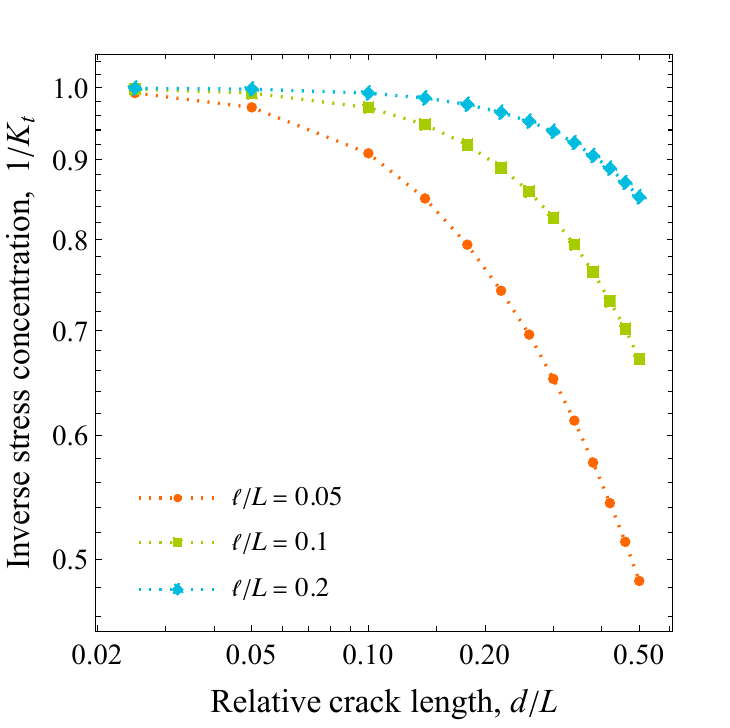}
	\caption{Dependence of inverse stress concentration at the crack tip on the relative length of crack under mode I (a, $K_t = \tau_{22}/t$) and mode II (b, $K_t = \tau_{12}/t$) loading conditions}
	\label{fig8}
\end{figure}

Alternatively, the failure analysis of pre-cracked bodies within SGE can be performed based on the phenomenological consideration of regularized solution for Cauchy stresses as suggested in Refs. \cite{askes2015understanding,vasiliev2019estimation,vasiliev2021failure}. For presented mode I and mode II problems with central crack, the influence of crack size on the inverse values of maximum Cauchy stress concentration is presented in Fig.\ref{fig8}. The stress concentration is evaluated at the crack tip for maximum values of normal stress $K_t = \tau_{22}/t$ in the mode I problem and for shear stress $K_t = \tau_{12}/t$ in the mode II problem. The presented values of inverse stress concentration can be treated as an assessment on the normalized nominal strength of pre-cracked bodies. For small cracks, which size is comparable to the length scale parameter $\ell$, SGE solution predicts the increase of normalized nominal strength up to unit value that corresponds to the absence of stress concentration. For the long cracks, the increase of stress concentration arises and its rate depends on the length scale parameter. Thus, the identification of $\ell$ can be also performed considering the dependence of the specimens' nominal strength on the crack size and specimen size \cite{vasiliev2021failure,vasiliev2019estimation,askes2011gradient}.


\section{Conclusions}
\label{con}

Conventional fifth-order approximation (that was used previously within SGE) implies the continuity of displacements and their first and second gradients. Such approximation does provide accurate numerical solutions within SGE, though for the crack problems one needs to use a very dense meshes to capture the discontinuity in the second gradient of displacements around the crack tip. Suggested  enriched finite elements embodies the appropriate SGE asymptotic solutions for cracks and allows us to obtained the accurate results on the coarse meshes. Simultaneously, we can evaluate the amplitude factors and related values of J-integral that can be used to predict the fracture of pre-cracked bodies accounting for the size effects within SGE. 

Proposed scheme \eqref{ue}-\eqref{Nii} allows us to preserve $C^1$-continuity of solution between the enriched and standard elements. This scheme is new and has not been discussed previously. Presented examples of calculations validates the efficiency of presented method. Moreover, the presented data on the dependence of amplitude factors of SGE asymptotic solutions on the crack length is the first systematic result in the field (for the best of authors knowledge). 

In the future, wide class of problems and related data on size effects for amplitude factors, J-integral and concentration of regularized stresses can be obtained following developed method. These studies would be an important step for further identification of material constants and validation of SGE  theory. The more general examples of problems with mixed mode loading and 3D formulation should be also considered.

%

\section*{Appendix A. Asymptotic solution}
\label{appA}
\setcounter{equation}{0}\renewcommand\theequation{A.\arabic{equation}}

The asymptotic solution used in the present study has been derived within the simplified SGE in Refs. \cite{gourgiotis2009plane,aravas2009plane}. Recently, we have re-examined this solution within the analysis of the high-order terms \cite{solyaev2024higher} and found that one can explicitly extract the classical part of this solution that correspond to the classical term $\sim r^{3/2}$ in Williams series. With such explicit definition of classical part, the displacement asymptotic solution for the mode I crack problem was derived in polar coordinates in the following form:

\begin{equation}
\label{A1}
\begin{aligned}
	\textbf u_I &= \textbf u^{(c)} +  \textbf u^{(g)}
\end{aligned}
\end{equation}
where the classical part $\textbf u^{(c)} = u_r^{(c)} \textbf e_r + u_\theta^{(c)} \textbf e_\theta$ is given by:

\begin{equation}
\label{A2}
\begin{aligned}
	u_r^{(c)} &= \tfrac{K_1}{4\mu} r^{3/2} 
	\left(
	(2\eta-3) \cos \tfrac{\theta}{2}
	+ \cos \tfrac{5\theta}{2}
	\right),\\
	u_\theta^{(c)} &= \tfrac{K_1}{4\mu} r^{3/2} 
	\left(
	(2\eta-3) \sin \tfrac{\theta}{2}
	+ \sin \tfrac{5\theta}{2}
	\right)
\end{aligned}
\end{equation}
and the additional gradient part $\textbf u^{(g)} = u_r^{(g)} \textbf e_r + u_\theta^{(g)} \textbf e_\theta$ is the following:
\begin{equation}
\label{A3}
\begin{aligned}
	u_r^{(g)} &= \tfrac{K_2}{4\mu} r^{3/2} 
	\left(
	\tfrac{4\eta-1}{2} \cos \tfrac{3\theta}{2}
	- \tfrac{8\eta+17}{6}\cos \tfrac{5\theta}{2}
	\right),\\
	u_\theta^{(g)} &= \tfrac{K_2}{4\mu} r^{3/2} 
	\left(
	-\tfrac{4\eta+1}{2} \sin \tfrac{3\theta}{2}
	+ \tfrac{8\eta+17}{6}\sin \tfrac{5\theta}{2}
	\right)
\end{aligned}
\end{equation}

Here $\textbf u^{(c)}$ \eqref{A2} is a standard term from Williams series with asymptotic behaviour $\sim r^{3/2}$ \cite{Williams1960}. In our derivations in Ref. \cite{solyaev2024higher}, the amplitude factors $K_1$ and $K_2$  in \eqref{A2}, \eqref{A3} were defined as $K_3$ and $K_{-1}$, respectively (due the ordering of these terms in the generalized Williams series within SGE). The form of asymptotic solution given in Ref. \cite{gourgiotis2009plane} can be obtain if we put $K_1 = A_1$ and $K_2 = 6(A_1-A_2)/(17+8\eta)$ ($A_1, A_2$ are the amplitude factors used in Ref. \cite{gourgiotis2009plane}).

For the mode II crack problem, the displacement solution in polar coordinates can be presented in the similar decomposed form:
\begin{equation}
\label{A4}
\begin{aligned}
	\textbf u_{II} &= \textbf u^{(c)} +  \textbf u^{(g)}
\end{aligned}
\end{equation}
with the classical part $\textbf u^{(c)} = u_r^{(c)} \textbf e_r + u_\theta^{(c)} \textbf e_\theta$ given by:
\begin{equation}
\label{A5}
\begin{aligned}
	u_r^{(c)} &= \tfrac{K_3}{4\mu} r^{3/2} 
	\left(
	(3-2\eta) \sin \tfrac{\theta}{2}
	- 5 \sin \tfrac{5\theta}{2}
	\right),\\
	u_\theta^{(c)} &= \tfrac{K_3}{4\mu} r^{3/2} 
	\left(
	(3+2\eta) \cos \tfrac{\theta}{2}
	- 5 \cos \tfrac{5\theta}{2}
	\right)
\end{aligned}
\end{equation}
and with the gradient part $\textbf u^{(g)} = u_r^{(g)} \textbf e_r + u_\theta^{(g)} \textbf e_\theta$ defined by:
\begin{equation}
\label{A6}
\begin{aligned}
	u_r^{(g)} &= \tfrac{K_4}{4\mu} r^{3/2} \sin \tfrac{\theta}{2}
	+\tfrac{K_3}{4\mu} r^{3/2}\left(
	\tfrac{3(1-4\eta)}{2}\sin \tfrac{3\theta}{2}
	+ \tfrac{23+8\eta}{2}\sin \tfrac{5\theta}{2}
	\right),\\
	u_\theta^{(g)} &= -\tfrac{K_4}{4\mu} r^{3/2} \cos \tfrac{\theta}{2}
	-\tfrac{K_3}{4\mu} r^{3/2}\left(
	\tfrac{3(1+4\eta)}{2}\cos \tfrac{3\theta}{2}
	- \tfrac{23+8\eta}{2}\cos \tfrac{5\theta}{2}
	\right)
\end{aligned}
\end{equation}

In contrast to the mode I solution \eqref{A1}-\eqref{A3}, the amplitudes of the mode II problem persist both in the gradient part \eqref{A6}. Nevertheless, only one of them ($K_3$) defines the classical term from the Williams series \eqref{A5}. The reason for the coupling of amplitudes in the gradient part of the mode II solution is an open question for the further research. Possibly, there exists another kind of simplified gradient theory (not in Aifantis form), which constitutive equations allow the total separation of the classical and the gradient parts of solution for the mode II crack problems. The form of asymptotic solution given in Ref. \cite{gourgiotis2009plane}, can be obtained from \eqref{A5}-\eqref{A6} if we put: $K_3 = 2B_2/(13+8\eta)$, $K_4 = B_1 - 2B_2(3-2\eta)/(13+8\eta)$ (where $B_1, B_2$ are the amplitude factors used in Ref. \cite{gourgiotis2009plane} for the mode II problem).  

The final form of the asymptotic solution $\textbf u = u_x\, \textbf e_x + u_y\, \textbf e_y$ \eqref{as} used for the enrichment of elements' shape functions can be obtained by using combination of solutions \eqref{A1}-\eqref{A6} and standard transformation between the polar $(r,\theta)$ and Cartesian coordinates $(x,y)$. Components of displacement vector in Cartesian coordinates in the main text of this paper were denoted as $u_x=u$, $u_y = v$ (see \eqref{as}).

\section*{Appendix B. Shape functions}
\label{appB}
\setcounter{equation}{0}\renewcommand\theequation{B.\arabic{equation}}

The compact representation of the shape functions for the Bell triangle can be given in terms of areal coordinates in the following form \cite{dasgupta1990higher}:
\begin{equation}
\label{Ni}
\begin{aligned}
	N_1 & = L_1^5 + 5L_1^4L_2 + 5L_1^4L_3 + 10L_1^3L_2^2 +10L_1^3L_3^2\\
		& +20L_1^3L_2L_3 + 30 r_{21}L_1^2L_2L_3^2 + 30 r_{31}L_1^2L_3L_2^2\\
	N_2 & = c_3L_1^4L_2 - c_2L_1^4L_3 + 4c_3L_1^3L_2^2 - 4c_2L_1^3L_3^2\\
		& + 4(c_3-c_2)L_1^3L_2L_3 - (3c_1+15r_{21}c_2)L_1^2L_2L_3^2\\
		& + (3c_1+15r_{31}c_3)L_1^2L_3L_2^2\\
	N_3 & = -b_3L_1^4L_2 + b_2L_1^4L_3 - 4b_3L_1^3L_2^2 + 4b_2L_1^3L_3^2\\
		& + 4(b_2-b_3)L_1^3L_2L_3 + (3b_1+15r_{21}b_2)L_1^2L_2L_3^2\\
		&-(3b_1+15r_{31}b_3)L_1^2L_3L_2^2\\
	N_4 & = \frac{c_3^2}{2}L_1^3L_2^2 + \frac{c_2^2}{2}L_1^3L_3^2 - c_2c_3L_1^3L_2L_3\\
		& + (c_1c_2 + \tfrac{5}{2}r_{21}c_2^2)L_2L_3^2L_1^2 + (c_1c_3 + \tfrac{5}{2}r_{31}c_3^2)L_3L_2^2L_1^2\\
	N_5 & = -b_3c_3L_1^3L_2^2 - b_2c_2L_1^3L_3^2 + (b_2c_3+b_3c_2)L_1^3L_2L_3\\
		& - (b_1c_2 + b_2c_1 + 5r_{21}b_2c_2)L_2L_3^2L_1^2\\
		& - (b_1c_3 + b_3c_1 + 5r_{31}b_3c_3)L_3L_2^2L_1^2\\
	N_6 & = \frac{b_3^2}{2}L_1^3L_2^2 + \frac{b_2^2}{2}L_1^3L_3^2 - b_2b_3L_1^3L_2L_3\\
		& +(b_1b_2 + \tfrac{5}{2}r_{21}b_2^2)L_2L_3^2L_1^2
		+(b_1b_3 + \tfrac{5}{2}r_{31}b_3^2)L_3L_2^2L_1^2
\end{aligned}
\end{equation}
where $L_i$ ($i=1,2,3$) are the areal coordinates of an element that are related to the global Cartesian coordinates $(x,y)$ as follows 

\begin{equation}
\label{Li}
\begin{aligned}
	x &= L_i x_i, \quad y = L_i y_i, \quad L_1 + L_2 + L_3 = 1\\
	L_i &= \frac{1}{\Delta} (a_i + b_ix + c_i y), \quad(i=1,2,3)\\
	\Delta &= (x_2-x_1)(y_3-y_1)-(x_3-x_1)(y_2-y_1)
\end{aligned}
\end{equation}
where $(x_i,y_i)$ ($i=1,2,3$) being the global Cartesian coordinates of the nodes of triangular element and 
\begin{equation}
\label{ai}
\begin{aligned}
	a_i = x_j y_k- x_ky_j, \quad b_i = y_j-y_k, \quad c_i = x_k-x_j
\end{aligned}
\end{equation}
where $i, j, k$ being cyclic permutations of 1, 2 and 3.

The presented six shape functions \eqref{Ni} define the degrees of freedom (displacements and their first and second derivatives) in the first node $(x_1,y_1)$ of the element. The remaining twelve shape functions $N_7...N_{18}$ correspond to the degrees of freedom at the second and the third nodes and they can be obtained by the cyclic permutations of the subscripts 1, 2, 3 in relations \eqref{Ni}.

The derivatives of functions given in terms of areal coordinates can be found by using the following relations (with summation over repeated indices):
\begin{equation}
\label{dL}
\begin{aligned}
	\frac{\partial(...)}{\partial x} = \frac{b_i}{\Delta}  \frac{\partial(...)}{\partial L_i},
	\qquad
	\frac{\partial(...)}{\partial y} = \frac{c_i}{\Delta} \frac{\partial(...)}{\partial L_i}
\end{aligned}
\end{equation}

\section*{References}
\bibliography{refs.bib}

\end{document}